\documentclass[reqno,12pt,a4paper]{amsart}

\usepackage{mathrsfs}
\usepackage{amssymb}
\usepackage{amsfonts}
\usepackage{latexsym}
\usepackage{amsthm}
\usepackage{graphicx}

\def\lb{\label}

\newcommand{\er}[1]{\textrm{(\ref{#1})}}

\begin{document}


\renewcommand{\theequation}{\arabic{section}.\arabic{equation}}
\theoremstyle{plain}
\newtheorem{theorem}{\bf Theorem}[section]
\newtheorem{lemma}[theorem]{\bf Lemma}
\newtheorem{corollary}[theorem]{\bf Corollary}
\newtheorem{proposition}[theorem]{\bf Proposition}
\newtheorem{definition}[theorem]{\bf Definition}
\newtheorem{example}[theorem]{\bf Example}

\theoremstyle{remark}
\newtheorem*{remarks}{\bf Remarks}
\newtheorem*{remark}{\bf Remark}

\def\a{\alpha}  \def\cA{{\mathcal A}}     \def\bA{{\bf A}}  \def\mA{{\mathscr A}}
\def\b{\beta}   \def\cB{{\mathcal B}}     \def\bB{{\bf B}}  \def\mB{{\mathscr B}}
\def\g{\gamma}  \def\cC{{\mathcal C}}     \def\bC{{\bf C}}  \def\mC{{\mathscr C}}
\def\G{\Gamma}  \def\cD{{\mathcal D}}     \def\bD{{\bf D}}  \def\mD{{\mathscr D}}
\def\d{\delta}  \def\cE{{\mathcal E}}     \def\bE{{\bf E}}  \def\mE{{\mathscr E}}
\def\D{\Delta}  \def\cF{{\mathcal F}}     \def\bF{{\bf F}}  \def\mF{{\mathscr F}}
\def\c{\chi}    \def\cG{{\mathcal G}}     \def\bG{{\bf G}}  \def\mG{{\mathscr G}}
\def\z{\zeta}   \def\cH{{\mathcal H}}     \def\bH{{\bf H}}  \def\mH{{\mathscr H}}
\def\e{\eta}    \def\cI{{\mathcal I}}     \def\bI{{\bf I}}  \def\mI{{\mathscr I}}
\def\p{\psi}    \def\cJ{{\mathcal J}}     \def\bJ{{\bf J}}  \def\mJ{{\mathscr J}}
\def\vT{\Theta} \def\cK{{\mathcal K}}     \def\bK{{\bf K}}  \def\mK{{\mathscr K}}
\def\k{\kappa}  \def\cL{{\mathcal L}}     \def\bL{{\bf L}}  \def\mL{{\mathscr L}}
\def\l{\lambda} \def\cM{{\mathcal M}}     \def\bM{{\bf M}}  \def\mM{{\mathscr M}}
\def\L{\Lambda} \def\cN{{\mathcal N}}     \def\bN{{\bf N}}  \def\mN{{\mathscr N}}
\def\m{\mu}     \def\cO{{\mathcal O}}     \def\bO{{\bf O}}  \def\mO{{\mathscr O}}
\def\n{\nu}     \def\cP{{\mathcal P}}     \def\bP{{\bf P}}  \def\mP{{\mathscr P}}
\def\r{\rho}    \def\cQ{{\mathcal Q}}     \def\bQ{{\bf Q}}  \def\mQ{{\mathscr Q}}
\def\s{\sigma}  \def\cR{{\mathcal R}}     \def\bR{{\bf R}}  \def\mR{{\mathscr R}}
\def\S{\Sigma}  \def\cS{{\mathcal S}}     \def\bS{{\bf S}}  \def\mS{{\mathscr S}}
\def\t{\tau}    \def\cT{{\mathcal T}}     \def\bT{{\bf T}}  \def\mT{{\mathscr T}}
\def\f{\phi}    \def\cU{{\mathcal U}}     \def\bU{{\bf U}}  \def\mU{{\mathscr U}}
\def\F{\Phi}    \def\cV{{\mathcal V}}     \def\bV{{\bf V}}  \def\mV{{\mathscr V}}
\def\P{\Psi}    \def\cW{{\mathcal W}}     \def\bW{{\bf W}}  \def\mW{{\mathscr W}}
\def\o{\omega}  \def\cX{{\mathcal X}}     \def\bX{{\bf X}}  \def\mX{{\mathscr X}}
\def\x{\xi}     \def\cY{{\mathcal Y}}     \def\bY{{\bf Y}}  \def\mY{{\mathscr Y}}
\def\X{\Xi}     \def\cZ{{\mathcal Z}}     \def\bZ{{\bf Z}}  \def\mZ{{\mathscr Z}}

\def\vr{\varrho}
\def\be{{\bf e}} \def\bc{{\bf c}}
\def\bx{{\bf x}} \def\by{{\bf y}}
\def\bv{{\bf v}} \def\bu{{\bf u}}
\def\Om{\Omega} \def\bp{{\bf p}}
\def\bbD{\pmb \Delta}
\def\mm{\mathrm m}
\def\mn{\mathrm n}

\newcommand{\mc}{\mathscr {c}}

\newcommand{\gA}{\mathfrak{A}}          \newcommand{\ga}{\mathfrak{a}}
\newcommand{\gB}{\mathfrak{B}}          \newcommand{\gb}{\mathfrak{b}}
\newcommand{\gC}{\mathfrak{C}}          \newcommand{\gc}{\mathfrak{c}}
\newcommand{\gD}{\mathfrak{D}}          \newcommand{\gd}{\mathfrak{d}}
\newcommand{\gE}{\mathfrak{E}}
\newcommand{\gF}{\mathfrak{F}}           \newcommand{\gf}{\mathfrak{f}}
\newcommand{\gG}{\mathfrak{G}}           
\newcommand{\gH}{\mathfrak{H}}           \newcommand{\gh}{\mathfrak{h}}
\newcommand{\gI}{\mathfrak{I}}           \newcommand{\gi}{\mathfrak{i}}
\newcommand{\gJ}{\mathfrak{J}}           \newcommand{\gj}{\mathfrak{j}}
\newcommand{\gK}{\mathfrak{K}}            \newcommand{\gk}{\mathfrak{k}}
\newcommand{\gL}{\mathfrak{L}}            \newcommand{\gl}{\mathfrak{l}}
\newcommand{\gM}{\mathfrak{M}}            \newcommand{\gm}{\mathfrak{m}}
\newcommand{\gN}{\mathfrak{N}}            \newcommand{\gn}{\mathfrak{n}}
\newcommand{\gO}{\mathfrak{O}}
\newcommand{\gP}{\mathfrak{P}}             \newcommand{\gp}{\mathfrak{p}}
\newcommand{\gQ}{\mathfrak{Q}}             \newcommand{\gq}{\mathfrak{q}}
\newcommand{\gR}{\mathfrak{R}}             \newcommand{\gr}{\mathfrak{r}}
\newcommand{\gS}{\mathfrak{S}}              \newcommand{\gs}{\mathfrak{s}}
\newcommand{\gT}{\mathfrak{T}}             \newcommand{\gt}{\mathfrak{t}}
\newcommand{\gU}{\mathfrak{U}}             \newcommand{\gu}{\mathfrak{u}}
\newcommand{\gV}{\mathfrak{V}}             \newcommand{\gv}{\mathfrak{v}}
\newcommand{\gW}{\mathfrak{W}}             \newcommand{\gw}{\mathfrak{w}}
\newcommand{\gX}{\mathfrak{X}}               \newcommand{\gx}{\mathfrak{x}}
\newcommand{\gY}{\mathfrak{Y}}              \newcommand{\gy}{\mathfrak{y}}
\newcommand{\gZ}{\mathfrak{Z}}             \newcommand{\gz}{\mathfrak{z}}

\def\ve{\varepsilon}   \def\vt{\vartheta}    \def\vp{\varphi}    \def\vk{\varkappa}

\def\A{{\mathbb A}} \def\B{{\mathbb B}} \def\C{{\mathbb C}}
\def\dD{{\mathbb D}} \def\E{{\mathbb E}} \def\dF{{\mathbb F}} \def\dG{{\mathbb G}}
\def\H{{\mathbb H}}\def\I{{\mathbb I}} \def\J{{\mathbb J}} \def\K{{\mathbb K}}
 \def\dL{{\mathbb L}}\def\M{{\mathbb M}} \def\N{{\mathbb N}} \def\O{{\mathbb O}}
 \def\dP{{\mathbb P}} \def\R{{\mathbb R}}\def\S{{\mathbb S}} \def\T{{\mathbb T}}
  \def\U{{\mathbb U}} \def\V{{\mathbb V}}\def\W{{\mathbb W}} \def\X{{\mathbb X}}
   \def\Y{{\mathbb Y}} \def\Z{{\mathbb Z}}


\def\la{\leftarrow}              \def\ra{\rightarrow}            \def\Ra{\Rightarrow}
\def\ua{\uparrow}                \def\da{\downarrow}
\def\lra{\leftrightarrow}        \def\Lra{\Leftrightarrow}


\def\lt{\biggl}                  \def\rt{\biggr}
\def\ol{\overline}               \def\wt{\widetilde}
\def\ul{\underline}
\def\no{\noindent}


\let\ge\geqslant                 \let\le\leqslant
\def\lan{\langle}                \def\ran{\rangle}
\def\/{\over}                    \def\iy{\infty}
\def\sm{\setminus}               \def\es{\emptyset}
\def\ss{\subset}                 \def\ts{\times}
\def\pa{\partial}                \def\os{\oplus}
\def\om{\ominus}                 \def\ev{\equiv}
\def\iint{\int\!\!\!\int}        \def\iintt{\mathop{\int\!\!\int\!\!\dots\!\!\int}\limits}
\def\el2{\ell^{\,2}}             \def\1{1\!\!1}
\def\sh{\sharp}
\def\wh{\widehat}
\def\bs{\backslash}
\def\intl{\int\limits}

\def\na{\mathop{\mathrm{\nabla}}\nolimits}
\def\sh{\mathop{\mathrm{sh}}\nolimits}
\def\ch{\mathop{\mathrm{ch}}\nolimits}
\def\where{\mathop{\mathrm{where}}\nolimits}
\def\all{\mathop{\mathrm{all}}\nolimits}
\def\as{\mathop{\mathrm{as}}\nolimits}
\def\Area{\mathop{\mathrm{Area}}\nolimits}
\def\arg{\mathop{\mathrm{arg}}\nolimits}
\def\const{\mathop{\mathrm{const}}\nolimits}
\def\det{\mathop{\mathrm{det}}\nolimits}
\def\diag{\mathop{\mathrm{diag}}\nolimits}
\def\diam{\mathop{\mathrm{diam}}\nolimits}
\def\dim{\mathop{\mathrm{dim}}\nolimits}
\def\dist{\mathop{\mathrm{dist}}\nolimits}
\def\Im{\mathop{\mathrm{Im}}\nolimits}
\def\Iso{\mathop{\mathrm{Iso}}\nolimits}
\def\Ker{\mathop{\mathrm{Ker}}\nolimits}
\def\Lip{\mathop{\mathrm{Lip}}\nolimits}
\def\rank{\mathop{\mathrm{rank}}\limits}
\def\Ran{\mathop{\mathrm{Ran}}\nolimits}
\def\Re{\mathop{\mathrm{Re}}\nolimits}
\def\Res{\mathop{\mathrm{Res}}\nolimits}
\def\res{\mathop{\mathrm{res}}\limits}
\def\sign{\mathop{\mathrm{sign}}\nolimits}
\def\span{\mathop{\mathrm{span}}\nolimits}
\def\supp{\mathop{\mathrm{supp}}\nolimits}
\def\Tr{\mathop{\mathrm{Tr}}\nolimits}
\def\BBox{\hspace{1mm}\vrule height6pt width5.5pt depth0pt \hspace{6pt}}


\newcommand\nh[2]{\widehat{#1}\vphantom{#1}^{(#2)}}
\def\dia{\diamond}

\def\Oplus{\bigoplus\nolimits}



\def\qqq{\qquad}
\def\qq{\quad}
\let\ge\geqslant
\let\le\leqslant
\let\geq\geqslant
\let\leq\leqslant
\newcommand{\ca}{\begin{cases}}
\newcommand{\ac}{\end{cases}}
\newcommand{\ma}{\begin{pmatrix}}
\newcommand{\am}{\end{pmatrix}}
\renewcommand{\[}{\begin{equation}}
\renewcommand{\]}{\end{equation}}
\def\eq{\begin{equation}}
\def\qe{\end{equation}}
\def\[{\begin{equation}}
\def\bu{\bullet}

\title[Two-sided estimates of  total bandwidth for
Schr\"odinger operators ]{Two-sided estimates of total bandwidth
for Schr\"odinger operators on periodic graphs}

\date{\today}
\author[Evgeny Korotyaev]{Evgeny Korotyaev}
\address{Saint-Petersburg State University, Universitetskaya nab. 7/9, St. Petersburg,
 199034, Russia,
and HSE University, 3A Kantemirovskaya ulitsa, St. Petersburg, 194100,
Russia, \ korotyaev@gmail.com, \
e.korotyaev@spbu.ru,}
\author[Natalia Saburova]{Natalia Saburova}
\address{Northern (Arctic) Federal University, Severnaya Dvina Emb. 17, Arkhangelsk,
 163002, Russia,
 \ n.saburova@gmail.com, \ n.saburova@narfu.ru}


\begin{abstract}
We consider Schr\"odinger operators with periodic potentials on periodic discrete graphs. Their spectrum consists of a finite number of bands. We obtain two-sided estimates of the total bandwidth for the Schr\"odinger operators in terms of geometric parameters of the graph and the potentials. In particular, we show that these estimates are sharp. It means that these estimates become identities for specific graphs and potentials. The proof is based on the Floquet theory and trace formulas for fiber operators. The traces are expressed as finite Fourier series of the quasimomentum with coefficients depending on the potentials and cycles of the quotient graph from some specific cycle sets. In order to obtain our results we estimate these Fourier coefficients in terms of geometric parameters of the graph and the potentials.
\end{abstract}

\maketitle

\section {\lb{Sec1}Introduction}
\setcounter{equation}{0}
Laplace and Schr\"odinger operators on graphs have a lot of
applications in physics,  chemistry and engineering. They are used
to study properties of different periodic media, e.g. nanomedia, see \cite{NG04} (two authors of
this survey, Novoselov and Geim, won the Nobel Prize for
discovering graphene). We consider Schr\"odinger operators $H=-\D+V$
on periodic graphs, where $\D$ is  the Laplacian   and $V$ is a
periodic potential. On a discrete graph, the Laplace operator acts
in the space of functions defined on the set of vertices of
the graph. Here vertices of the graph play a crucial role, while
edges show the interaction between the vertices.  Basically, two
types of the discrete Laplace operator have been studied: the
combinatorial  Laplacian  and the normalized Laplacian. The spectrum
of $H$ is a finite  union  of spectral bands $\s_j$, $j=1,\ldots,\n$.
Our main goal is to obtain two-sided estimates for the total bandwidth
$\sum_{j=1}^\n|\s_j|$ in terms of geometric parameters of the graphs and the
potentials.

There are a lot of results about spectral properties for the
one-dimensional case. Most of the results were obtained for the lattice $\Z$, since there are applications to Toda lattices, see, e.g., \cite{T89}. Here corresponding operators of the Lax pairs are Jacobi matrices with periodic coefficients. We mention that Last \cite{L92} determined two-sided sharp estimates \er{1diS} for the total bandwidth of Schr\"odinger operators with periodic potentials on $\Z$, see also \cite{DS83}, \cite{KK03} and references therein.

In the continuous case we know only estimates for Schr\"odinger
operators with periodic potentials on the real line. Korotyaev
determined two-sided estimates for gap lengths of Schr\"odinger
operators with periodic potentials \cite{K98}, \cite{K03}. In
\cite{K00} Korotyaev obtained various estimates of spectral bands
and, in particular, estimates of action variables for the KdV
equation on the circle. The case of matrix-valued potentials was
discussed by Chelkak and Korotyaev \cite{CK06}. It is important that
the proof for these cases was based on  the conformal mapping theory
(associated with quasimomentum) and trace formulas.

We discuss the multidimensional case. We do not know estimates for
the continuous case.  We consider  the discrete case, where there
are estimates of the total bandwidth. It is known that the Lebesgue
measure $|\s(\D)|$ of the spectrum of the combinatorial Laplacian
$\D$ on periodic graphs satisfies $|\s(\D)|\leq2\vk_+$, where
$\vk_+$ is the maximum vertex degree of the graph. But, in contrast
to the one-dimensional case, the bands may overlap. Thus, the total
bandwidth may exceed this number. \emph{Upper} estimates of the
total bandwidth for the Schr\"odinger operator with a periodic
potential in terms of geometric parameters of the graph were
obtained in \cite{KS14,KS18,KS20}. For Schr\"odinger operators
with periodic magnetic potentials similar upper estimates were
obtained in \cite{KS17}. We do not know results about \emph{lower}
estimates of the total bandwidth of the Laplace and Schr\"odinger
operators on periodic graphs.

In this paper we consider Laplace and Schr\"odinger operators with
periodic  potentials on periodic discrete graphs. We obtain
two-sided estimates of the total bandwidth for the Schr\"odinger
operators in terms of geometric parameters of the graph and the
potentials. We show that these estimates are sharp, i.e., they
become identities for specific graphs and potentials. The proof is
based on the Floquet theory and recent trace formulas for fiber
operators from \cite{KS21}.

\subsection{Schr\"odinger operators on periodic graphs}
Let $\cG=(\cV,\cE)$ be a connected infinite graph, possibly  having
loops and multiple edges and embedded into the space $\R^d$. Here
$\cV$  is the set of its vertices and $\cE$ is the set of its
unoriented edges. Considering each edge in $\cE$ to have two
orientations, we introduce the set $\cA$ of all oriented edges. An
edge starting at a vertex $x$ and ending at a vertex $y$ from $\cV$
will be denoted as the ordered pair $(x,y)\in\cA$ and is said to be
\emph{incident} to the vertices. Let $\ul\be=(y,x)$ be the inverse
edge of $\be=(x,y)\in\cA$. Vertices $x,y\in\cV$ will be called
\emph{adjacent} and denoted by $x\sim y$, if $(x,y)\in \cA$. We
define the degree $\vk_x$ of the vertex $x\in\cV$ as the number of
all edges from $\cA$, starting at $x$.

Let $\G$ be a lattice of rank $d$ in $\R^d$ with a basis $\{a_1,\ldots,a_d\}$, i.e.,
$$
\G=\Big\{a : a=\sum_{s=1}^dn_sa_s, \; (n_s)_{s=1}^d\in\Z^d\Big\},
$$
and define the \emph{fundamental cell} $\Omega$ of the lattice $\G$
by
\[\lb{fuce}
\Omega=\Big\{\bx\in\R^d : \bx=\sum_{s=1}^d\bx_sa_s, \; (\bx_s)_{s=1}^d\in
[0,1)^d\Big\}.
\]
We introduce the equivalence relation on $\R^d$:
$$
\bx\equiv \by \; (\hspace{-4mm}\mod \G) \qq\Leftrightarrow\qq \bx-\by\in\G \qqq
\forall\, \bx,\by\in\R^d.
$$

We consider \emph{locally finite $\G$-periodic graphs} $\cG$, i.e., graphs satisfying the
following conditions:
\begin{itemize}
  \item[1)] $\cG=\cG+a$ for any $a\in\G$;
  \item[2)] the quotient graph  $\cG_*=\cG/\G$ is finite.
\end{itemize}
The basis $a_1,\ldots,a_d$ of the lattice $\G$ is called the {\it
periods}  of $\cG$. We also call the quotient graph $\cG_*=\cG/\G$
the \emph{fundamental graph} of the periodic graph $\cG$. The
fundamental graph $\cG_*$ is a graph on the $d$-dimensional torus
$\R^d/\G$. The graph $\cG_*=(\cV_*,\cE_*)$ has the vertex set
$\cV_*=\cV/\G$, the set $\cE_*=\cE/\G$ of unoriented edges and the
set $\cA_*=\cA/\G$ of oriented edges which are finite.

Let $\ell^2(\cV)$ be the Hilbert space of all square summable
functions  $f:\cV\to \C$ equipped with the norm
$$
\|f\|^2_{\ell^2(\cV)}=\sum_{x\in\cV}|f_x|^2<\infty.
$$
We consider the Schr\"odinger operator $H$ acting on the Hilbert space $\ell^2(\cV)$ and given by
\[
\lb{Sh} H=H_o+V,\qqq H_o=-\D,
\]
where the potential $V$ is real valued and satisfies for all
$(x,a)\in\cV\ts\G$:
\[
\lb{Pot} (V f)_x=V_xf_x, \qqq V_{x+a}=V_x,
\]
and $\D$ is the \emph{combinatorial Laplacian} on
$f\in\ell^2(\cV)$ defined by
\[\lb{DLO}
(\D f)_x=\sum_{(x,y)\in\cA}\big(f_x-f_y\big), \qqq x\in\cV.
\]
The sum in \er{DLO} is taken over all oriented edges starting at the
vertex $x$. It is well known   that $\D$ is self-adjoint and the
point 0 belongs to its spectrum $\s(\D)$ contained in $[0,2\vk_+]$,
i.e.,
\[
\lb{bf}
0\in\s(\D)\subseteq[0,2\vk_+],\qqq
\textrm{where}\qqq
\vk_+=\max_{x\in\cV}\vk_x<\iy,
\]
and $\vk_x$ is the degree of the vertex $x$, see, e.g., \cite{MW89}.

\subsection{Edge and cycle indices} Here we will define an {\it edge index}.
It was introduced in \cite{KS14} in order to express fiber
Laplacians and Schr\"odinger operators on the fundamental graph in
terms of edge indices, see \er{fLao}. For each $\bx\in\R^d$ we
introduce the vector $\bx_\A\in\R^d$ by
\[
\lb{cola}
\bx_\A=(\bx_1,\ldots,\bx_d), \qqq \textrm{where} \qq \bx=\textstyle\sum\limits_{s=1}^d\bx_s a_s,
\]
i.e., $\bx_\A$ is the coordinate vector of $\bx$ with respect to the
basis  $\A=\{a_j\}_{j=1}^d$ of the lattice~$\G$. For any vertex
$x\in\cV$ of a $\G$-periodic graph $\cG$ the following unique
representation holds true:
\[\lb{Dv}
x=x_0+[x], \qq \textrm{where}\qq (x_0,[x])\in\Omega\ts\G,
\]
and $\Omega$ is a fundamental cell of the lattice $\G$ defined by
\er{fuce}. In other words, each vertex $x$ can be obtained from a
vertex  $x_0\in \Omega$ by a shift by a vector $[x]\in\G$. For any
oriented edge $\be=(x,y)\in\cA$ we define the \emph{edge index}
$\t(\be)$ as the vector of the lattice $\Z^d$ given by
\[
\lb{in}
\t(\be)=[y]_\A-[x]_\A\in\Z^d,
\]
where $[\cdot]\in\G$ is defined by \er{Dv} and the vector
$[x]_\A\in\Z^d$ is given by \er{cola}.

On the set $\cA$ of all oriented edges of the $\G$-periodic graph $\cG$ we define the surjection
\[\lb{sur}
\gf:\cA\rightarrow\cA_*=\cA/\G,
\]
which maps each $\be\in\cA$ to its equivalence class
$\be_*=\gf(\be)$ which is  an oriented edge of the fundamental graph
$\cG_*$. For any oriented edge $\be_*\in\cA_*$  we define the
\emph{edge index} $\t(\be_*)\in\Z^d$ by
\[
\lb{dco}
\t(\be_*)=\t(\be) \qq \textrm{ for some $\be\in\cA$ \; such that }
  \; \be_*=\gf(\be), \qqq \be_*\in\cA_*,
\]
where $\gf$ is defined by \er{sur}. Thus, edge indices of the
fundamental graph $\cG_*$ are induced by edge indices of the
periodic graph~$\cG$. The edge index $\t(\be_*)$ is uniquely
determined by \er{dco}, since
$$
\t(\be+a)=\t(\be),\qqq \forall\, (\be,a)\in\cA \ts\G.
$$

A \emph{path} $\bp$ in a graph $\cG=(\cV,\cA)$ is a sequence of
consecutive edges
\[\lb{depa}
\bp=(\be_1,\be_2,\ldots,\be_n), \qqq \where \qq
\be_s=(x_{s-1},x_s)\in\cA,\qq s=1,\ldots,n,
\]
for some vertices $x_0,x_1,\ldots,x_n\in\cV$. The vertices $x_0$ and
$x_n$ are called the \emph{initial} and \emph{terminal} vertices of
the path $\bp$, respectively. If $x_0=x_n$, then the path $\bp$ is
called a \emph{cycle}. The number $n$ of edges in a cycle $\bc$ is
called the \emph{length} of $\bc$ and is denoted by $|\bc|$, i.e.,
$|\bc|=n$. The \emph{reverse} of the path $\bp$ given by \er{depa}
is the path $\ul\bp=(\ul\be_n,\ldots,\ul\be_1)$.

\begin{remark} A path $\bp$ is uniquely defined by the sequence of
its oriented edges $(\be_1,\be_2,\ldots,\be_n)$. The sequence of
its vertices $(x_0,x_1,\ldots,x_n)$ does not uniquely define~$\bp$,
since multiple edges are allowed in the graph $\cG$.
\end{remark}

Let $\cC$ be the set of all cycles of the fundamental graph
$\cG_*$.  For any cycle $\bc\in\cC$ we define the \emph{cycle index}
$\t(\bc)\in\Z^d$ by
\[\lb{cyin}
\t(\bc)=\sum\limits_{\be\in\bc}\t(\be),  \qqq \bc\in\cC.
\]
From the definition of indices it follows that
\[\lb{inin}
\t(\ul\be\,)=-\t(\be), \qq \forall\,\be\in\cA_* \qqq \textrm{and}\qqq
\t(\ul\bc\,)=-\t(\bc), \qq \forall\,\bc\in\cC.
\]

\begin{figure}[h]
\centering
\unitlength 1.0mm 
\linethickness{0.4pt}
\ifx\plotpoint\undefined\newsavebox{\plotpoint}\fi 

\begin{picture}(50,43)(0,0)
\hspace{-10mm}
\put(30,2){$\Downarrow$}
\bezier{20}(24.5,16.3)(24.5,22.3)(24.5,28.3)
\bezier{20}(25.0,16.6)(25.0,22.6)(25.0,28.6)
\bezier{20}(25.5,16.9)(25.5,22.9)(25.5,28.9)
\bezier{20}(26.0,17.2)(26.0,23.2)(26.0,29.2)
\bezier{20}(26.5,17.5)(26.5,23.5)(26.5,29.5)
\bezier{20}(27.0,17.8)(27.0,23.8)(27.0,29.8)
\bezier{20}(27.5,18.1)(27.5,24.1)(27.5,30.1)
\bezier{20}(28.0,18.4)(28.0,24.4)(28.0,30.4)
\bezier{20}(28.5,18.7)(28.5,24.7)(28.5,30.7)
\bezier{20}(29.0,19.0)(29.0,25.0)(29.0,31.0)

\bezier{20}(29.5,19.3)(29.5,25.3)(29.5,31.3)
\bezier{20}(30.0,19.6)(30.0,25.6)(30.0,31.6)
\bezier{20}(30.5,19.9)(30.5,25.9)(30.5,31.9)
\bezier{20}(31.0,20.2)(31.0,26.2)(31.0,32.2)
\bezier{20}(31.5,20.5)(31.5,26.5)(31.5,32.5)
\bezier{20}(32.0,20.8)(32.0,26.8)(32.0,32.8)
\bezier{20}(32.5,21.1)(32.5,27.1)(32.5,33.1)
\bezier{20}(33.0,21.4)(33.0,27.4)(33.0,33.4)
\bezier{20}(33.5,21.7)(33.5,27.7)(33.5,33.7)
\bezier{20}(34,22)(34,28)(34,34)

\put(17.5,29.5){$\scriptstyle v_2+a_2$}

\put(14,10){\circle{1}}
\put(28,10){\circle{1}}
\put(34,10){\circle{1}}
\put(48,10){\circle{1}}

\put(18,16){\circle{1}}
\put(24,16){\circle*{1}}
\put(38,16){\circle{1}}
\put(44,16){\circle{1}}

\put(14,22){\circle{1}}
\put(28,22){\circle*{1}}
\put(34,22){\circle{1}}
\put(48,22){\circle{1}}

\put(18,28){\circle{1}}
\put(24,28){\circle{1}}
\put(38,28){\circle{1}}
\put(44,28){\circle{1}}

\put(14,34){\circle{1}}
\put(28,34){\circle{1}}
\put(34,34){\circle{1}}
\put(48,34){\circle{1}}

\put(18,40){\circle{1}}
\put(24,40){\circle{1}}
\put(38,40){\circle{1}}
\put(44,40){\circle{1}}

\put(28,10){\line(1,0){6.00}}
\put(18,16){\line(1,0){6.00}}
\put(38,16){\line(1,0){6.00}}

\put(28,22){\line(1,0){6.00}}
\put(18,28){\line(1,0){6.00}}
\put(38,28){\line(1,0){6.00}}

\put(28,34){\line(1,0){6.00}}
\put(18,40){\line(1,0){6.00}}
\put(38,40){\line(1,0){6.00}}

\put(14,10){\line(2,3){4.00}}

\put(34,10){\line(2,3){4.00}}
\put(24,16){\line(2,3){4.00}}
\put(44,16){\line(2,3){4.00}}

\put(14,22){\line(2,3){4.00}}
\put(34,22){\line(2,3){4.00}}
\put(24,28){\line(2,3){4.00}}
\put(44,28){\line(2,3){4.00}}

\put(14,34){\line(2,3){4.00}}
\put(34,34){\line(2,3){4.00}}

\put(28,10){\line(-2,3){4.00}}
\put(48,10){\line(-2,3){4.00}}
\put(38,16){\line(-2,3){4.00}}
\put(18,16){\line(-2,3){4.00}}

\put(28,22){\line(-2,3){4.00}}
\put(48,22){\line(-2,3){4.00}}
\put(38,28){\line(-2,3){4.00}}
\put(18,28){\line(-2,3){4.00}}

\put(28,34){\line(-2,3){4.00}}
\put(48,34){\line(-2,3){4.00}}

\put(30,18){$\scriptstyle a_1$}
\put(20.5,22){$\scriptstyle a_2$}

\put(18,16){\line(1,0){6.00}}
\put(28,22){\line(1,0){6.00}}
\qbezier(28,22)(26,25)(24,28)
\qbezier(28,10)(26,13)(24,16)

\put(24,16){\vector(0,1){12.0}}
\put(33,21.3){\vector(3,2){0.5}}

\qbezier(24,16)(29,19)(34,22)
\put(30,28){$\Omega$}
\put(32,22){\vector(1,0){0.5}}
\put(26.4,19.6){\vector(2,3){0.5}}
\put(26,25){\vector(-2,3){0.5}}

\put(31,23){$\scriptstyle \be_2$}
\put(25.5,26){$\scriptstyle \be_3$}
\put(24.0,20.5){$\scriptstyle \be_1$}
\put(27.5,23){$\scriptstyle v_1$}
\put(35,21.5){$\scriptstyle v_2+a_1$}
\put(30,35){$\scriptstyle v_2+a_1+a_2$}

\put(20.5,13.8){$\scriptstyle O=v_2$}

\put(14,10){\line(-1,0){3.00}}
\put(14,22){\line(-1,0){3.00}}
\put(14,34){\line(-1,0){3.00}}
\put(48,10){\line(1,0){3.00}}
\put(48,22){\line(1,0){3.00}}
\put(48,34){\line(1,0){3.00}}

\bezier{50}(14,10)(15,8.5)(16,7)
\bezier{50}(34,10)(35,8.5)(36,7)
\bezier{50}(28,10)(27,8.5)(26,7)
\bezier{50}(48,10)(47,8.5)(46,7)

\bezier{50}(18,40)(17,41.5)(16,43)
\bezier{50}(38,40)(37,41.5)(36,43)
\bezier{50}(24,40)(25,41.5)(26,43)
\bezier{50}(44,40)(45,41.5)(46,43)
\put(7,38){$\bG$}
\put(55,26){$v_1,v_2\in\Omega$}
\end{picture}\hspace{5mm}
\begin{picture}(50,43)(0,0)
\put(30,2){$\Downarrow$}
\bezier{20}(26.5,19.3)(26.5,25.3)(26.5,31.3)
\bezier{20}(27.0,19.6)(27.0,25.6)(27.0,31.6)
\bezier{20}(27.5,19.9)(27.5,25.9)(27.5,31.9)
\bezier{20}(28.0,20.2)(28.0,26.2)(28.0,32.2)
\bezier{20}(28.5,20.5)(28.5,26.5)(28.5,32.5)
\bezier{20}(29.0,20.8)(29.0,26.8)(29.0,32.8)
\bezier{20}(29.5,21.1)(29.5,27.1)(29.5,33.1)
\bezier{20}(30.0,21.4)(30.0,27.4)(30.0,33.4)
\bezier{20}(30.5,21.7)(30.5,27.7)(30.5,33.7)
\bezier{20}(31.0,22.0)(31.0,28.0)(31.0,34.0)

\bezier{20}(31.5,22.3)(31.5,28.3)(31.5,34.3)
\bezier{20}(32.0,22.6)(32.0,28.6)(32.0,34.6)
\bezier{20}(32.5,22.9)(32.5,28.9)(32.5,34.9)
\bezier{20}(33.0,23.2)(33.0,29.2)(33.0,35.2)
\bezier{20}(33.5,23.5)(33.5,29.5)(33.5,35.5)
\bezier{20}(34.0,23.8)(34.0,29.8)(34.0,35.8)
\bezier{20}(34.5,24.1)(34.5,30.1)(34.5,36.1)
\bezier{20}(35.0,24.4)(35.0,30.4)(35.0,36.4)
\bezier{20}(35.5,24.7)(35.5,30.7)(35.5,36.7)
\bezier{20}(36,25)(36,31)(36,37)

\put(17.5,29.5){$\scriptstyle v_2-a_1$}

\put(14,10){\circle{1}}
\put(28,10){\circle{1}}
\put(34,10){\circle{1}}
\put(48,10){\circle{1}}

\put(18,16){\circle{1}}
\put(24,16){\circle{1}}
\put(38,16){\circle{1}}
\put(44,16){\circle{1}}

\put(14,22){\circle{1}}
\put(28,22){\circle*{1}}
\put(34,22){\circle{1}}
\put(48,22){\circle{1}}

\put(18,28){\circle{1}}
\put(24,28){\circle{1}}
\put(38,28){\circle{1}}
\put(44,28){\circle{1}}

\put(14,34){\circle{1}}
\put(28,34){\circle{1}}
\put(34,34){\circle*{1}}
\put(48,34){\circle{1}}

\put(18,40){\circle{1}}
\put(24,40){\circle{1}}
\put(38,40){\circle{1}}
\put(44,40){\circle{1}}

\put(28,10){\line(1,0){6.00}}
\put(18,16){\line(1,0){6.00}}
\put(38,16){\line(1,0){6.00}}

\put(28,22){\line(1,0){6.00}}
\put(18,28){\line(1,0){6.00}}
\put(38,28){\line(1,0){6.00}}

\put(28,34){\line(1,0){6.00}}
\put(18,40){\line(1,0){6.00}}
\put(38,40){\line(1,0){6.00}}

\put(14,10){\line(2,3){4.00}}

\put(34,10){\line(2,3){4.00}}
\put(24,16){\line(2,3){4.00}}
\put(44,16){\line(2,3){4.00}}

\put(14,22){\line(2,3){4.00}}
\put(34,22){\line(2,3){4.00}}
\put(24,28){\line(2,3){4.00}}
\put(44,28){\line(2,3){4.00}}

\put(14,34){\line(2,3){4.00}}
\put(34,34){\line(2,3){4.00}}

\put(28,10){\line(-2,3){4.00}}
\put(48,10){\line(-2,3){4.00}}
\put(38,16){\line(-2,3){4.00}}
\put(18,16){\line(-2,3){4.00}}

\put(28,22){\line(-2,3){4.00}}
\put(48,22){\line(-2,3){4.00}}
\put(38,28){\line(-2,3){4.00}}
\put(18,28){\line(-2,3){4.00}}

\put(28,34){\line(-2,3){4.00}}
\put(48,34){\line(-2,3){4.00}}

\put(28,34){\line(1,0){6.00}}
\put(28,22){\line(1,0){6.00}}
\qbezier(28,22)(26,25)(24,28)
\qbezier(38,40)(36,37)(34,34)
\qbezier(28,22)(26,19)(24,16)
\qbezier(34,34)(36,31)(38,28)

\put(33,22){\vector(1,0){0.5}}
\put(27.2,21.0){\vector(2,3){0.5}}
\put(24.4,27.4){\vector(-2,3){0.5}}
\put(31,20){$\scriptstyle \be_2$}
\put(23,25){$\scriptstyle \be_3$}
\put(23.0,19){$\scriptstyle \be_1$}
\put(25.5,17){$\scriptstyle O$}

\put(33.5,25.3){$\scriptstyle a_1$}
\put(26,28){$\scriptstyle a_2$}

\put(26,19){\vector(0,1){12.0}}
\put(35,24.3){\vector(3,2){0.5}}
\qbezier(26,19)(31,22)(36,25)

\put(30,27){$\Omega$}
\put(27,23.5){$\scriptstyle v_1$}
\put(35,21.5){$\scriptstyle v_2-a_2$}
\put(31.3,32){$\scriptstyle v_2$}

\put(19.8,14.0){$\scriptstyle v_2-a_1-a_2$}

\put(14,10){\line(-1,0){3.00}}
\put(14,22){\line(-1,0){3.00}}
\put(14,34){\line(-1,0){3.00}}
\put(48,10){\line(1,0){3.00}}
\put(48,22){\line(1,0){3.00}}
\put(48,34){\line(1,0){3.00}}

\bezier{50}(14,10)(15,8.5)(16,7)
\bezier{50}(34,10)(35,8.5)(36,7)
\bezier{50}(28,10)(27,8.5)(26,7)
\bezier{50}(48,10)(47,8.5)(46,7)

\bezier{50}(18,40)(17,41.5)(16,43)
\bezier{50}(38,40)(37,41.5)(36,43)
\bezier{50}(24,40)(25,41.5)(26,43)
\bezier{50}(44,40)(45,41.5)(46,43)

\put(7,38){$\bG$}
\end{picture}

\begin{picture}(30,30)(0,0)
\put(5,10){\circle*{1}}
\put(17,26){\circle*{1}}
\put(2,7.0){$v_2$}
\put(17,27){$v_1$}

\bezier{200}(5,10.0)(11,18.0)(17,26.0)
\bezier{200}(5,10)(2,27)(17,26)
\bezier{200}(5,10)(21,10)(17,26)
\put(9,30){$\bG_*$}

\put(18.5,15.5){$\be_2$}
\put(1.5,22){$\be_3$}
\put(8.0,19){$\be_1$}

\put(11.0,18.0){\vector(2,3){1.0}}
\put(7.0,23.0){\vector(-1,-1){1.0}}
\put(17.7,17.0){\vector(-1,-2){1.0}}
\put(16,12){$\scriptstyle(1,0)$}
\put(3,25.0){$\scriptstyle(0,1)$}
\put(9.5,15){$\scriptstyle(0,0)$}

\end{picture} \hspace{45mm}
\begin{picture}(20,30)(0,0)
\put(5,10){\circle*{1}}
\put(17,26){\circle*{1}}
\put(2,7.0){$v_2$}
\put(17,27){$v_1$}

\bezier{200}(5,10.0)(11,18.0)(17,26.0)
\bezier{200}(5,10)(2,27)(17,26)
\bezier{200}(5,10)(21,10)(17,26)
\put(9,30){$\bG_*$}

\put(18.5,15.5){$\be_2$}
\put(1.5,22){$\be_3$}
\put(8.0,19){$\be_1$}

\put(11.0,18.0){\vector(2,3){1.0}}
\put(7.0,23.0){\vector(-1,-1){1.0}}
\put(17.7,17.0){\vector(-1,-2){1.0}}
\put(16,12){$\scriptstyle(0,-1)$}
\put(3,25.5){$\scriptstyle(-1,0)$}
\put(9.5,15){$\scriptstyle(1,1)$}
\end{picture}
\vspace{-8mm}
\caption{\scriptsize The hexagonal lattice $\bG$ and its fundamental graph $\bG_*$ with edge indices; $a_1,a_2$ are the periods of $\bG$. The fundamental cell $\Omega$ is shaded. The vertices $v_1,v_2$ of $\bG$ from $\Omega$ are black points. Edge indices depend on the choice of the embedding of the periodic graph into $\R^2$ (i.e., the choice of $\Omega$). Cycle indices do not depend on this choice.} \lb{FGHei}
\end{figure}
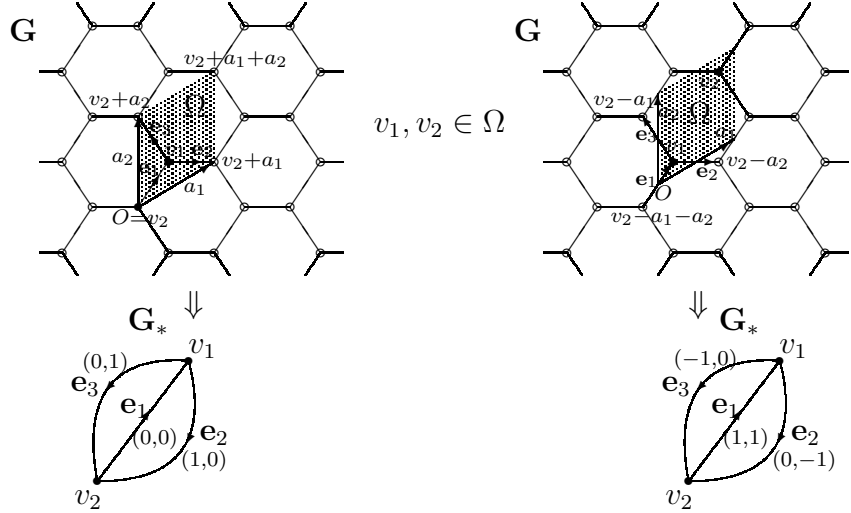

\begin{remarks} 1) Edge indices depend on the choice of the embedding
of the periodic graph $\cG$ into $\R^d$ (see Fig. \ref{FGHei}).
Cycle indices \emph{do not} depend on this choice. This property
of cycle indices has the following simple explanation. Any cycle
$\bc$ in the fundamental graph $\cG_*$ is obtained by factorization
of a path in the $\G$-periodic graph $\cG$ connecting some
$\G$-equivalent vertices $x\in\cV$ and $x+a\in\cV$, $a\in\G$. The
index of the cycle $\bc$ is equal to $\mm=(m_j)_{j=1}^d\in\Z^d$,
where $a=m_1a_1+\ldots+m_da_d$, and, therefore, does not depend on
the choice of the embedding.

2) The index of a cycle $\bc\in\cC$ represents the coordinates of $\bc$
with  respect to some basis of the space $\cC^+\ss\cC$ of cycles
with non-zero indices.
\end{remarks}

\subsection{Invariant $\cI$}
Edges with non-zero indices are called \emph{bridges}. They provide
connectivity of the periodic graph and removal of all bridges
disconnects the periodic graph into infinitely many connected
components (see Fig. \ref{FBri}). Let $\cB_*$ be the set of all
bridges of the fundamental graph $\cG_*=(\cV_*,\cA_*)$, i.e., $\cB_*=\supp\t$, where $\t:\cA_*\ra\Z^d$ is the index form defined by \er{in},
\er{dco}. The number of the fundamental graph bridges depends on the choice of the embedding of the periodic graph $\cG$ into the space $\R^d$, i.e., this number is not an invariant for $\cG$. We define the number
\[\lb{dIm}
\textstyle\cI=\frac12\,\min\limits_{\cG\ss\R^d}\#\cB_*, \qqq \cB_*=\supp\t,
\]
where the minimum is taken over all embeddings of the periodic graph
$\cG$ into $\R^d$, and $\#M$ denotes the number of elements in a set $M$. This number $\cI$ exists, since the fundamental
graph is finite and we described $\cI$  in Theorem 3.2 from
\cite{KS20}. For simple periodic graphs (the $d$-dimensional
lattice, the hexagonal lattice, the Kagome lattice, etc.) it is not
difficult to determine this number. But for an arbitrary periodic
graph this may be a rather complicated problem.

\setlength{\unitlength}{1.0mm}
\begin{figure}[h]
\centering
\unitlength 1.0mm 
\linethickness{0.4pt}
\ifx\plotpoint\undefined\newsavebox{\plotpoint}\fi 

\begin{picture}(50,50)(0,0)
\put(0,10.0){(\emph{a})}
\put(22.5,22.5){\vector(1,0){15.00}}
\put(22.5,22.5){\vector(0,1){15.00}}
\multiput(22.5,37.5)(4,0){4}{\line(1,0){2}}
\multiput(37.5,22.5)(0,4){4}{\line(0,1){2}}

\bezier{30}(22.5,22.5)(22.5,30)(22.5,37.5)
\bezier{30}(23.0,22.5)(23.0,30)(23.0,37.5)
\bezier{30}(23.5,22.5)(23.5,30)(23.5,37.5)
\bezier{30}(24.0,22.5)(24.0,30)(24.0,37.5)
\bezier{30}(24.5,22.5)(24.5,30)(24.5,37.5)
\bezier{30}(25.0,22.5)(25.0,30)(25.0,37.5)
\bezier{30}(25.5,22.5)(25.5,30)(25.5,37.5)
\bezier{30}(26.0,22.5)(26.0,30)(26.0,37.5)
\bezier{30}(26.5,22.5)(26.5,30)(26.5,37.5)
\bezier{30}(27.0,22.5)(27.0,30)(27.0,37.5)
\bezier{30}(27.5,22.5)(27.5,30)(27.5,37.5)
\bezier{30}(28.0,22.5)(28.0,30)(28.0,37.5)
\bezier{30}(28.5,22.5)(28.5,30)(28.5,37.5)
\bezier{30}(29.0,22.5)(29.0,30)(29.0,37.5)
\bezier{30}(29.5,22.5)(29.5,30)(29.5,37.5)
\bezier{30}(30.0,22.5)(30.0,30)(30.0,37.5)
\bezier{30}(30.5,22.5)(30.5,30)(30.5,37.5)
\bezier{30}(31.0,22.5)(31.0,30)(31.0,37.5)
\bezier{30}(31.5,22.5)(31.5,30)(31.5,37.5)
\bezier{30}(32.0,22.5)(32.0,30)(32.0,37.5)
\bezier{30}(32.5,22.5)(32.5,30)(32.5,37.5)
\bezier{30}(33.0,22.5)(33.0,30)(33.0,37.5)
\bezier{30}(33.5,22.5)(33.5,30)(33.5,37.5)
\bezier{30}(34.0,22.5)(34.0,30)(34.0,37.5)
\bezier{30}(34.5,22.5)(34.5,30)(34.5,37.5)
\bezier{30}(35.0,22.5)(35.0,30)(35.0,37.5)
\bezier{30}(35.5,22.5)(35.5,30)(35.5,37.5)
\bezier{30}(36.0,22.5)(36.0,30)(36.0,37.5)
\bezier{30}(36.5,22.5)(36.5,30)(36.5,37.5)
\bezier{30}(37.0,22.5)(37.0,30)(37.0,37.5)
\bezier{30}(37.5,22.5)(37.5,30)(37.5,37.5)

\put(34,20.0){$a_1$}
\put(18.5,35){$a_2$}
\put(28.5,29){$\Omega$}

\put(10,15){\line(1,-1){5.00}}
\put(10,15){\line(1,1){5.00}}
\put(20,15){\line(-1,-1){5.00}}
\put(20,15){\line(-1,1){5.00}}
\linethickness{1.5pt}
\put(20,15){\line(1,0){5.00}}
\put(15,20){\line(0,1){5.00}}
\put(10,15){\line(-1,0){2.00}}
\put(10,30){\line(-1,0){2.00}}
\put(10,45){\line(-1,0){2.00}}
\put(50,15){\line(1,0){2.00}}
\put(50,30){\line(1,0){2.00}}
\put(50,45){\line(1,0){2.00}}
\put(15,10){\line(0,-1){2.00}}
\put(30,10){\line(0,-1){2.00}}
\put(45,10){\line(0,-1){2.00}}
\put(15,50){\line(0,1){2.00}}
\put(30,50){\line(0,1){2.00}}
\put(45,50){\line(0,1){2.00}}

\linethickness{0.4pt}
\put(10,15){\circle{1.1}}
\put(15,10){\circle{1.1}}
\put(20,15){\circle{1.1}}
\put(15,20){\circle{1.1}}
\put(25,15){\line(1,-1){5.00}}
\put(25,15){\line(1,1){5.00}}
\put(35,15){\line(-1,-1){5.00}}
\put(35,15){\line(-1,1){5.00}}
\linethickness{1.5pt}
\put(35,15){\line(1,0){5.00}}
\put(30,20){\line(0,1){5.00}}
\linethickness{0.4pt}
\put(25,15){\circle{1.1}}
\put(30,10){\circle{1.1}}
\put(35,15){\circle{1.1}}
\put(30,20){\circle{1.1}}
\put(40,15){\line(1,-1){5.00}}
\put(40,15){\line(1,1){5.00}}
\put(50,15){\line(-1,-1){5.00}}
\put(50,15){\line(-1,1){5.00}}
\linethickness{1.5pt}
\put(45,20){\line(0,1){5.00}}
\linethickness{0.4pt}
\put(40,15){\circle{1.1}}
\put(45,10){\circle{1.1}}
\put(50,15){\circle{1.1}}
\put(45,20){\circle{1.1}}

\put(10,30){\line(1,-1){5.00}}
\put(10,30){\line(1,1){5.00}}
\put(20,30){\line(-1,-1){5.00}}
\put(20,30){\line(-1,1){5.00}}
\linethickness{1.5pt}
\put(20,30){\line(1,0){5.00}}
\put(15,35){\line(0,1){5.00}}
\linethickness{0.4pt}
\put(10,30){\circle{1.1}}
\put(15,25){\circle{1.1}}
\put(20,30){\circle{1.1}}
\put(15,35){\circle{1.1}}
\put(25,30){\line(1,-1){5.00}}
\put(25,30){\line(1,1){5.00}}
\put(35,30){\line(-1,-1){5.00}}
\put(35,30){\line(-1,1){5.00}}
\linethickness{1.5pt}
\put(35,30){\line(1,0){5.00}}
\put(30,35){\line(0,1){5.00}}
\linethickness{0.4pt}
\put(25,30){\circle*{1.5}}
\put(30,25){\circle*{1.5}}
\put(35,30){\circle*{1.5}}
\put(30,35){\circle*{1.5}}
\put(40,30){\line(1,-1){5.00}}
\put(40,30){\line(1,1){5.00}}
\put(50,30){\line(-1,-1){5.00}}
\put(50,30){\line(-1,1){5.00}}
\linethickness{1.5pt}
\put(45,35){\line(0,1){5.00}}
\linethickness{0.4pt}
\put(40,30){\circle{1.1}}
\put(45,25){\circle{1.1}}
\put(50,30){\circle{1.1}}
\put(45,35){\circle{1.1}}

\put(10,45){\line(1,-1){5.00}}
\put(10,45){\line(1,1){5.00}}
\put(20,45){\line(-1,-1){5.00}}
\put(20,45){\line(-1,1){5.00}}
\linethickness{1.5pt}
\put(20,45){\line(1,0){5.00}}
\linethickness{0.4pt}
\put(10,45){\circle{1.1}}
\put(15,40){\circle{1.1}}
\put(20,45){\circle{1.1}}
\put(15,50){\circle{1.1}}
\put(25,45){\line(1,-1){5.00}}
\put(25,45){\line(1,1){5.00}}
\put(35,45){\line(-1,-1){5.00}}
\put(35,45){\line(-1,1){5.00}}
\linethickness{1.5pt}
\put(35,45){\line(1,0){5.00}}
\linethickness{0.4pt}
\put(25,45){\circle{1.1}}
\put(30,40){\circle{1.1}}
\put(35,45){\circle{1.1}}
\put(30,50){\circle{1.1}}
\put(40,45){\line(1,-1){5.00}}
\put(40,45){\line(1,1){5.00}}
\put(50,45){\line(-1,-1){5.00}}
\put(50,45){\line(-1,1){5.00}}
\put(40,45){\circle{1.1}}
\put(45,40){\circle{1.1}}
\put(50,45){\circle{1.1}}
\put(45,50){\circle{1.1}}


\end{picture}\qqq
\begin{picture}(50,50)(0,0)
\put(0,10.0){(\emph{b})}
\put(22.5,22.5){\vector(1,0){15.00}}
\put(22.5,22.5){\vector(0,1){15.00}}
\multiput(22.5,37.5)(4,0){4}{\line(1,0){2}}
\multiput(37.5,22.5)(0,4){4}{\line(0,1){2}}

\put(34,20.0){$a_1$}
\put(18.5,35){$a_2$}
\put(28.5,29){$\Omega$}

\bezier{30}(22.5,22.5)(22.5,30)(22.5,37.5)
\bezier{30}(23.0,22.5)(23.0,30)(23.0,37.5)
\bezier{30}(23.5,22.5)(23.5,30)(23.5,37.5)
\bezier{30}(24.0,22.5)(24.0,30)(24.0,37.5)
\bezier{30}(24.5,22.5)(24.5,30)(24.5,37.5)
\bezier{30}(25.0,22.5)(25.0,30)(25.0,37.5)
\bezier{30}(25.5,22.5)(25.5,30)(25.5,37.5)
\bezier{30}(26.0,22.5)(26.0,30)(26.0,37.5)
\bezier{30}(26.5,22.5)(26.5,30)(26.5,37.5)
\bezier{30}(27.0,22.5)(27.0,30)(27.0,37.5)
\bezier{30}(27.5,22.5)(27.5,30)(27.5,37.5)
\bezier{30}(28.0,22.5)(28.0,30)(28.0,37.5)
\bezier{30}(28.5,22.5)(28.5,30)(28.5,37.5)
\bezier{30}(29.0,22.5)(29.0,30)(29.0,37.5)
\bezier{30}(29.5,22.5)(29.5,30)(29.5,37.5)
\bezier{30}(30.0,22.5)(30.0,30)(30.0,37.5)
\bezier{30}(30.5,22.5)(30.5,30)(30.5,37.5)
\bezier{30}(31.0,22.5)(31.0,30)(31.0,37.5)
\bezier{30}(31.5,22.5)(31.5,30)(31.5,37.5)
\bezier{30}(32.0,22.5)(32.0,30)(32.0,37.5)
\bezier{30}(32.5,22.5)(32.5,30)(32.5,37.5)
\bezier{30}(33.0,22.5)(33.0,30)(33.0,37.5)
\bezier{30}(33.5,22.5)(33.5,30)(33.5,37.5)
\bezier{30}(34.0,22.5)(34.0,30)(34.0,37.5)
\bezier{30}(34.5,22.5)(34.5,30)(34.5,37.5)
\bezier{30}(35.0,22.5)(35.0,30)(35.0,37.5)
\bezier{30}(35.5,22.5)(35.5,30)(35.5,37.5)
\bezier{30}(36.0,22.5)(36.0,30)(36.0,37.5)
\bezier{30}(36.5,22.5)(36.5,30)(36.5,37.5)
\bezier{30}(37.0,22.5)(37.0,30)(37.0,37.5)
\bezier{30}(37.5,22.5)(37.5,30)(37.5,37.5)
\put(10,15){\line(1,-1){5.00}}
\put(10,15){\line(1,1){5.00}}
\put(20,15){\line(-1,-1){5.00}}
\put(20,15){\line(-1,1){5.00}}

\linethickness{0.4pt}
\put(10,15){\circle{1.1}}
\put(15,10){\circle{1.1}}
\put(20,15){\circle{1.1}}
\put(15,20){\circle{1.1}}
\put(25,15){\line(1,-1){5.00}}
\put(25,15){\line(1,1){5.00}}
\put(35,15){\line(-1,-1){5.00}}
\put(35,15){\line(-1,1){5.00}}
\linethickness{0.4pt}
\put(25,15){\circle{1.1}}
\put(30,10){\circle{1.1}}
\put(35,15){\circle{1.1}}
\put(30,20){\circle{1.1}}
\put(40,15){\line(1,-1){5.00}}
\put(40,15){\line(1,1){5.00}}
\put(50,15){\line(-1,-1){5.00}}
\put(50,15){\line(-1,1){5.00}}
\linethickness{0.4pt}
\put(40,15){\circle{1.1}}
\put(45,10){\circle{1.1}}
\put(50,15){\circle{1.1}}
\put(45,20){\circle{1.1}}

\put(10,30){\line(1,-1){5.00}}
\put(10,30){\line(1,1){5.00}}
\put(20,30){\line(-1,-1){5.00}}
\put(20,30){\line(-1,1){5.00}}
\linethickness{0.4pt}
\put(10,30){\circle{1.1}}
\put(15,25){\circle{1.1}}
\put(20,30){\circle{1.1}}
\put(15,35){\circle{1.1}}
\put(25,30){\line(1,-1){5.00}}
\put(25,30){\line(1,1){5.00}}
\put(35,30){\line(-1,-1){5.00}}
\put(35,30){\line(-1,1){5.00}}
\linethickness{0.4pt}
\put(25,30){\circle*{1.1}}
\put(30,25){\circle*{1.1}}
\put(35,30){\circle*{1.1}}
\put(30,35){\circle*{1.1}}
\put(40,30){\line(1,-1){5.00}}
\put(40,30){\line(1,1){5.00}}
\put(50,30){\line(-1,-1){5.00}}
\put(50,30){\line(-1,1){5.00}}
\linethickness{0.4pt}
\put(40,30){\circle{1.1}}
\put(45,25){\circle{1.1}}
\put(50,30){\circle{1.1}}
\put(45,35){\circle{1.1}}

\put(10,45){\line(1,-1){5.00}}
\put(10,45){\line(1,1){5.00}}
\put(20,45){\line(-1,-1){5.00}}
\put(20,45){\line(-1,1){5.00}}
\linethickness{0.4pt}
\put(10,45){\circle{1.1}}
\put(15,40){\circle{1.1}}
\put(20,45){\circle{1.1}}
\put(15,50){\circle{1.1}}
\put(25,45){\line(1,-1){5.00}}
\put(25,45){\line(1,1){5.00}}
\put(35,45){\line(-1,-1){5.00}}
\put(35,45){\line(-1,1){5.00}}
\linethickness{0.4pt}
\put(25,45){\circle{1.1}}
\put(30,40){\circle{1.1}}
\put(35,45){\circle{1.1}}
\put(30,50){\circle{1.1}}
\put(40,45){\line(1,-1){5.00}}
\put(40,45){\line(1,1){5.00}}
\put(50,45){\line(-1,-1){5.00}}
\put(50,45){\line(-1,1){5.00}}
\put(40,45){\circle{1.1}}
\put(45,40){\circle{1.1}}
\put(50,45){\circle{1.1}}
\put(45,50){\circle{1.1}}

\end{picture}
\vspace{-0.7cm} \caption{ \footnotesize  \emph{a}) A periodic graph $\cG$; $a_1,a_2$ are the periods of $\cG$; the fundamental cell $\Omega$ is shaded; bold edges are bridges of $\cG$; \emph{b}) removal of all bridges
disconnects the graph $\cG$ into infinitely many connected components.}
\lb{FBri}
\end{figure}
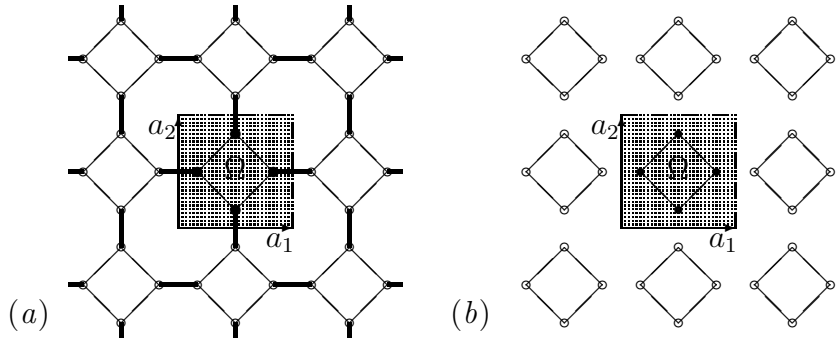

We recall some properties of the number $\cI$ (see Theorem 2.1.\emph{v} and Proposition 2.2 in \cite{KS20}).

\begin{proposition}\lb{InPr}
i) The number $\cI$ is an \textbf{invariant} of the $\G$-periodic graph $\cG$, i.e., it does not depend on the choice of

$\bu$ the embedding of $\cG$ into $\R^d$;

$\bu$ the basis $a_1,\ldots,a_d$ of the lattice $\G$.

ii) The invariant $\cI$ satisfies
\[\lb{cIsa}
d\leq\cI\leq\b, \qqq \b=\#\cE_*-\#\cV_*+1,
\]
where $d$ is the rank of the lattice $\G$ and $\b$ is the
Betti number of the fundamental graph $\cG_*=(\cV_*,\cE_*)$.

Moreover, for any $n\in\N$ there exists a periodic graph such that $\b-\cI=n$.
\end{proposition}

\subsection{Spectrum of Schr\"odinger operators.} We introduce the Hilbert space
\[\lb{Hisp}
\mH=L^2\Big(\T^{d},{dk\/(2\pi)^d}\,,\ell^2(\cV_*)\Big)=\int_{\T^{d}}^{\os}\ell^2(\cV_*)\,{dk
\/(2\pi)^d}\,, \qqq \T^d=\R^d/(2\pi\Z)^d,
\]
i.e., a constant fiber direct integral, equipped with the norm
$
\|g\|^2_{\mH}=\int_{\T^d}\|g(k)\|_{\ell^2(\cV_*)}^2\frac{dk}{(2\pi)^d}\,,
$
where the function $g(k)\in\ell^2(\cV_*)$ for almost all
$k\in\T^d$. Here
$\ell^2(\cV_*)=\C^\nu$ is the fiber space, $\nu=\#\cV_*$. The parameter
$k\in\T^d$ is called the \emph{quasimomentum}. We recall Theorem 1.1 from \cite{KS14}.

\begin{theorem}\label{TFDC}
The combinatorial Schr\"odinger operator $H=-\D+V$ on $\ell^2(\cV)$
has  the following decomposition into a constant fiber direct
integral
$$
UH U^{-1}=\int^\oplus_{\T^d}H(k)\,{dk\/(2\pi)^d}\,,
$$
for some unitary operator $U:\ell^2(\cV)\to\mH$ (the Gelfand
transform).  The fiber Schr\"odinger operator $H(k)$ on
$\ell^2(\cV_*)$ is given by
\[
\label{Hk'}
H(k)=H_o(k)+V, \qqq H_o(k)=-\D(k),\qqq \forall\,k\in\T^d.
\]
Here $V$ is the potential on $\ell^2(\cV_*)$, and $\D(k)$ is the fiber Laplacian having the form
\[
\label{fLao}
\big(\D(k)f\big)_x=\sum_{\be=(x,y)\in\cA_*}\big(f_x-e^{i\lan\t(\be),\,k\ran}f_y\big),
\qqq f\in\ell^2(\cV_*),\qqq x\in \cV_*,
\]
where $\t(\be)$ is the index of the edge $\be\in\cA_*$ defined by
\er{in},  \er{dco}, and $\lan\cdot\,,\cdot\,\ran$ denotes the
standard inner product in $\R^d$.
\end{theorem}

\begin{remarks} 1) The fiber Laplacian is expressed in terms of edge
indices which are not invariant and depend on the choice of the
embedding of the periodic graph $\cG$ into $\R^d$. But the invariance of indices of the fundamental graph cycles yields that the fiber Laplacians \er{fLao} written with respect to distinct embeddings of $\cG$ into $\R^d$ are unitarily equivalent under some gauge transform. Nevertheless, the proper embedding simplifies analysis of the operator $H$.

2) We can consider the fiber Laplacian $\D(k)$, $k\in\T^{d}$,  as a
discrete magnetic Laplacian with the magnetic vector potential
$\a(\be)=\lan\t(\be),k\ran$, $\be\in\cA_*$, on the fundamental
graph $\cG_*$.

3) The minimal number of exponents $e^{i\lan\t(\be),\,\cdot\,\ran}\neq1$, $\be\in\cA_*$, in the identities \er{fLao} for the fiber Laplacians $\D(\cdot)$ is equal to $2\cI$, where the invariant $\cI$ is defined by \er{dIm}.
\end{remarks}

Each fiber operator $H(k)$, $k\in\T^{d}$, has $\n=\#\cV_*$ real eigenvalues
$\l_j(k)$, $j\in\N_\n=\{1,\ldots,\n\}$, which are labeled in
non-decreasing order (counting multiplicities) by
$$
\l_{1}(k)\leq\l_{2}(k)\leq\ldots\leq\l_{\nu}(k),
\qqq \forall\,k\in\T^{d}.
$$
Each $\l_j(\cdot)$, $j\in\N_\n$, is a real and piecewise analytic
function  on the torus $\T^{d}$ and creates the \emph{spectral band}
$\s_j(H)$ given by
\[\lb{ban.1H}
\s_j(H)=[\l_j^-,\l_j^+]=\l_j(\T^{d}).
\]
Note that $\l_\n^+=\l_\n(0)$ (see \cite{SS92}). If
$\l_j(\cdot)=\L_j=\const$  on some subset of $\T^d$ of positive
Lebesgue measure, then the operator $H$ on $\cG$ has the eigenvalue
$\L_j$ of infinite multiplicity. We call $\{\L_j\}$ a \emph{flat
band}. Thus, the spectrum of the Schr\"odinger operator $H$ on the
periodic graph $\cG$ has the form
$$
\s(H)=\bigcup_{k\in\T^d}\s\big(H(k)\big)=
\bigcup_{j=1}^{\nu}\s_j(H)=\s_{ac}(H)\cup \s_{fb}(H),
$$
where $\s_{ac}(H)$ is the absolutely continuous spectrum, which is a
union of non-degenerate bands from \er{ban.1H}, and $\s_{fb}(H)$ is the set of
all flat bands (eigenvalues of infinite multiplicity).

We remark  that the last spectral band $\s_\n(H)$ of the
Schr\"odinger operator $H=-\D+V$ on periodic graphs is
non-degenerate, see \cite{KS19}. Moreover, if $\n\geq2$ and there
are no loops with non-zero index in the fundamental graph $\cG_*$,
then at least two spectral bands of $H$ are non-degenerate (see
Proposition 2.5 in \cite{KS14}).

\medskip

The paper is organized as follows. In Section \ref{Sec2} we formulate our main results:

$\bu$ estimates of the total bandwidth for the combinatorial Laplace and Schr\"odinger operators with periodic potentials on periodic graphs (Theorem \ref{TVECS} and Corollary \ref{CVECS});

$\bu$ estimates of the total bandwidth for the normalized Laplace operators (Theorem \ref{TVENL}).

Section \ref{Sec3} is devoted to trace formulas for Laplace and Schr\"odinger operators on periodic graphs which are used in the proof of the main results. In Section \ref{Sec4} we prove the estimates of the total bandwidth for the combinatorial Laplace and Schr\"odinger operators. In Section \ref{Sec5} we estimate the total bandwidth for the normalized Laplace operators. In Section \ref{Sec6} we apply the obtained results and estimate the total bandwidth for Laplacians on some specific periodic graphs.

\section{Main results}
\setcounter{equation}{0} \lb{Sec2}

\subsection{Overview of the total bandwidth estimates}
In order to present our results we describe known total bandwidth
estimates for Schr\"odinger operators $H=H_o+V$ with periodic
potentials $V$, where the unperturbed operator $H_o=-\D$. There are
only \emph{upper} estimates determined in
\cite{KS14,KS17,KS18,KS20}. We do not know lower estimates.

Define the total bandwidth $\gS(H)$ of an operator $H$ by
$$
\gS(H)= \sum_{j=1}^{\n}|\s_j(H)|.
$$

In \cite{KS14} the authors estimated the total bandwidth $\gS(H)$ of the Schr\"odinger operator $H$ in terms of geometric parameters of the graph:
\[
\lb{esLm0} \gS(H)\le 4b, \qqq \textrm{where} \qqq
b=\textstyle\frac12\,\#\supp\t,
\]
and $\t:\cA_*\ra\Z^d$ is the index form defined by \er{in},
\er{dco}. This number $b$ depends essentially on the choice of the
embedding of the periodic graph $\cG$ into the space $\R^d$, i.e.,
$b$ is not an invariant for $\cG$.

Similar estimates for the normalized Laplacian $\D_\gn$ (see \er{DNLA}) were obtained in \cite{KS19}:
\[\lb{esNL}
\gS (\D_\gn) \leq 2\gb,\qqq \gb=\sum_{x\in
\cV_*}\frac{\gb_x}{\vk_x}\,,
\]
where $\gb_x$ is the number of the fundamental graph bridges (i.e., edges with non-zero indices) starting at the vertex $x\in\cV_*$, and $\vk_x$ is the
degree of $x$.

For the Schr\"odinger operators $H$ we derived the following
estimate in  \cite{KS17}:
\[
\lb{esLm1} \gS(H)\le 4\b,
\]
where $\b=\#\cE_*-\#\cV_*+1$ is the Betti number of the fundamental graph $\cG_*=(\cV_*,\cE_*)$. Similar estimates were obtained for magnetic Schr\"odinger operators with periodic magnetic potentials in \cite{KS17}.

Finally in \cite{KS20} the authors introduced the new invariant
$\cI$, see \er{dIm}, and proved that
\[
\lb{ues1} \gS(H)\leq 4\cI.
\]

Let us compare the estimates \er{esLm0}, \er{esLm1} and \er{ues1}.
Due to \er{dIm} and \er{cIsa}, the invariant $\cI$ satisfies
$$
\cI=\min\limits_{\cG\ss\R^d} b\qqq\textrm{and}\qqq \cI\leq\b.
$$
Moreover, the difference between the Betti number $\b$ and the
invariant $\cI$ may be any non-negative integer number (for specific graphs), see Proposition \ref{InPr}.\emph{ii}). The number $b$ can be less, equal or significantly greater (dependently on the embedding of the periodic graph into $\R^d$) than the Betti number $\b$. Thus, \er{ues1} gives the
best upper estimate. Remark that \er{ues1} is sharp since we have an identity for specific graphs and potentials, and all estimates do not depend on potentials $V$.

In \cite{KS19} it was shown that the last band
$\s_\n(H)=[\l_\n^-,\l_\n^+]$  of the Schr\"odinger operator
$H=H_o+V$ is non-degenerate, i.e., $\l_\n^-<\l_\n^+$, and the
following estimate holds:
$$
0<C^{-1}|\s_\n(H_o)|\leq|\s_\n(H)|\leq
C\,|\s_\n(H_o)|,
$$
for a constant $C>0$ depending on geometric parameters of the graph
and the potential $V$.

Estimates on effective masses for Laplace  and Schr\"odinger
operators on periodic graphs were determined in a series of papers
\cite{K08}, \cite{KS16}, \cite{KS19}. A band localization for Laplace and Schr\"odinger operators on periodic graphs was discussed in \cite{FLP17}, \cite{KS15}, \cite{LP08}.

\subsection{Estimates for combinatorial Laplace and Schr\"odinger operators}
In order to describe our main results we need the following definitions. Let
$w=(w_x)_{x\in \cV_*}$ be a real potential on $\cV_*$. Define the
minimum and maximum and the diameter of the potential $w$ by
\[
\lb{vl+-} w_-=\min_{x\in\cV_*}w_x,\qqq
w_+=\max_{x\in\cV_*}w_x,\qqq \diam w=w_+- w_-.
\]

We present two-sided estimates of the total bandwidth for the Schr\"odinger
operators $H$. Our main goal is to estimate the total bandwidth for $H$ from \emph{below}.

\begin{theorem}\lb{TVECS}
Let $H=-\D+V$ be the Schr\"odinger operator defined by \er{Sh} --
\er{DLO}  on a periodic graph $\cG$. Then its total bandwidth
$\gS(H)= \sum_{j=1}^{\n}|\s_j(H)|$  satisfies
\[
\lb{CEC1} \frac{2d_*}{v_*^{\n-1}}\leq \gS(H)\leq 4\cI,
\]
where the invariant $\cI$ is defined by \er{dIm},
\[\lb{CEC11}
v_*=\vk_++\diam(V-\vk), \qqq d_*=\left\{\begin{array}{cl}
              d, & \textrm{if $d$ is even} \\
              d+1, & \textrm{if $d$ is odd}
            \end{array}\right.,
\]
$\vk_+$ and $\diam(V-\vk)$ are given by \er{vl+-}. In particular,
the Lebesgue measure $|\s(H)|$ of the spectrum of $H$ satisfies
\[
\lb{LMH}
\frac{2d_*}{\n v_*^{\n-1}}\leq|\s(H)|\leq 4\cI.
\]

\end{theorem}

\begin{remarks}
1) The lower estimate in \er{CEC1} gives another
proof of the known fact that the absolutely continuous spectrum of
$H$ is not empty.

2) The upper estimate in \er{CEC1} was proved in \cite{KS20}. The
proof of the lower estimate in \er{CEC1} consists of the following
steps:

\no $\bu$ The fiber Schr\"odinger operators are represented as fiber
adjacency operators on a modified  quotient graph with additional
loops and  with specific weights.

\no $\bu$ In Corollary \ref{CHEs}, using the spectral theorem, we
estimate the total bandwidth for  the Schr\"odinger operator $H$ in
terms of the total bandwidth of $H^n$ for any $n\in\N$.

\no $\bu$ In Theorem \ref{TEsS}, using the trace formulas for the
$n$-th power of  the fiber Schr\"odinger operator obtained in
\cite{KS21}, we show that the total bandwidth of $H^n$ is more than $2\cN_n^{odd}$, where $\cN_n^{odd}$ is the number of all cycles of length $n$ with an odd sum of index components in the fundamental graph $\cG_*$.

\no $\bu$ In Lemma \ref{LFtZ}, we show that the cycle indices form
the lattice $\Z^d$. Using properties of cycle indices, we show that
the number $\cN_n^{odd}\geq nd_*$ for some optimal $n\leq\n$, where
$d_*$ is given in \er{CEC11}.

3) Trace formulas for the discrete Laplacians on (mostly regular)
graphs and the numbers of graph cycles were analyzed by various
authors, see \cite{Ah87,Br91,CJK15,Mn07,TW03} and references
therein.
\end{remarks}

Taking $V=0$ in \er{CEC1}, we obtain similar estimates
for the unperturbed operator $H_o=-\D$. The eigenvalues of the
unperturbed fiber operator $H_o(k)$ will be denoted by $\l^o_j(k)$,
$j\in\N_\n$. The spectral bands $\s_j(H_o)$, $j\in\N_{\n}$, for the
operator $H_o$ have the form $\s_j(H_o)=[\l_j^{o-},\l_j^{o+}]=\l_j^o(\T^d)$.

\begin{corollary}\lb{CVECS}
Let $H_o=-\D$, where $\D$ is the combinatorial Laplacian defined by
\er{DLO} on a periodic graph  $\cG$. Then its total bandwidth
$\gS(H_o)$  satisfies
\[
\lb{CED1}
\frac{2d_*}{\vk_*^{\n-1}}\leq\gS(H_o)\leq 4\cI,\qqq
\textrm{where}\qqq \qq \vk_*=2\vk_+-\vk_-,
\]
and $\vk_\pm$ have the form \er{vl+-}.
\end{corollary}

\begin{remark}
If $\cG$ is a regular graph of degree $\vk_+$, i.e., all vertices of
$\cG$ have the same degree $\vk_+$, then $\vk_*=\vk_+$. For non-regular graphs $\vk_*>\vk_+$. Thus, among all periodic graphs with the same $d_*$, $\vk_+$ and $\n$, the regular one has the largest lower bound on the total bandwidth $\gS(H_o)$ in \er{CED1}. This agrees with the fact that deleting an edge from the fundamental graph $\cG_*$ increases the lower endpoint $\l_j^{o-}$ of each spectral band $\s_j(H_o)$, $j\in\N_\n$, of the operator $H_o=-\D$ (in particular, the last spectral band $\s_\n(H_o)=[\l_\n^{o-},0]$ shrinks), see, e.g. Corollary 4.2 in \cite{FLP20}.
\end{remark}

In order to show that the estimates \er{CEC1} are sharp we consider
Schr\"odinger operators $H=-\D+V$ with periodic potentials $V$ on
the simplest periodic graph, the one-dimensional lattice $\Z$. We
estimate the total bandwidth of $H$ using Theorem \ref{TVECS}. The
following corollary shows that in this case the estimates \er{CEC1}, \er{CEC11} coincide with the known Last's estimates \cite{L92}
for Schr\"odinger operators on $\Z$.

\begin{corollary}\lb{ESO}
Consider the Schr\"odinger operator $H=-\D+V$  on $\Z$, where $\D$
is the Laplacian on $\Z$ given by
$$
(-\D f)_x=f_{x+1}+f_{x-1}-2f_x, \qqq f\in\ell^2(\Z), \qq x\in\Z,
$$
and the potential $V$ is real $\n$-periodic, $V_{x+\n}=V_x$,
$x\in\Z$, $\n\geq3$. Then
\[\lb{1diS}
\frac{4}{v_*^{\n-1}}\leq\big|\s(H)\big|=\gS(H)\leq4, \qqq
\textrm{where} \qqq v_*=2+\diam V.
\]
\end{corollary}

\begin{remarks}
1) Last's proof of the lower estimate in \er{1diS} (see \cite{L92}) is based on the theory of periodic scalar Jacobi operators. For such operators the relation between the energy $\l$ and the one-dimensional quasimomentum $k$ (the dispersion relation) has a very specific form:
$$
D(\l)=2\cos k,
$$
where $D(\l)$ is the so called discriminant (or Lyapunov function). For an arbitrary periodic graph the dispersion relation may have a very complicated form and is difficult to analyze. Our proof of the lower estimates is based on the trace formulas for the Schr\"odinger operators on periodic graphs and bounds on the number of the fundamental graph cycles from some specific sets.

2) Consider the Schr\"odinger operator $H_t=-\D+tV$ from Corollary
\ref{ESO}, where the coupling constant $t>1$ and in addition
$V_x\neq V_y$ for all distinct $x,y\in \N_\n$. Then each band
$\s_j(H_t)$, $j\in\N_\n$, of the operator $H_t=-\D+tV$   has the
following asymptotics
$$
\big|\s_j(H_t)\big|=\frac4{t^{\n-1}}\prod_{x\in \N_\n, x\neq
j}|V_x-V_j|^{-1}+{O(1)\/t^\n}\,, \qqq \textrm{as}\qq t\to\iy,
$$
see \cite{KKu04}. Thus, the estimate \er{1diS} and the lower
estimate in \er{CEC1} are sharp.

3) Consider the Schr\"odinger  operators $H_t=-\D+tV$ with periodic
potentials $V$ on a periodic graph $\cG$, where the coupling
constant $t>1$ and in addition $V_x\neq V_y$ for all distinct
$x,y\in \cV_*$. Then  the following asymptotics holds true (see
Theorem 2.2 in \cite{KS14}):
$$
|\s(H_t)|=C+O(1/t),\qqq \textrm{as}\qq t\to\iy,
$$
and $C=0$ iff there are no loops with non-zero index in the fundamental graph $\cG_*$.
\end{remarks}

\subsection{Estimates for normalized Laplacians}

We define a \emph{normalized Laplacian} $\D_\gn$ on $\ell^2(\cV)$ by
\begin{equation}\lb{DNLA}
(\D_\gn f)_x=\frac1{\sqrt{\vk_x}}\sum_{(x,y)\in\cA}
\bigg(\frac{f_x}{\sqrt{\vk_x}}-\frac{f_y}{\sqrt{\vk_y}}\bigg)\,,
\qqq f\in\ell^2(\cV),\qqq x\in\cV.
\end{equation}
It is known (see, e.g., \cite{MW89}) that the normalized Laplacian
$\D_\gn$ is a bounded self-adjoint operator on $\ell^2(\cV)$ and its
spectrum $\s(\D_\gn)$ is a subset of the segment $[0,2]$, containing
the point 0, i.e.,
$$
0\in\s(\D_\gn)\subseteq[0,2].
$$

\begin{remark}
If $\cG$ is a regular graph of degree $\vk_+$, then the normalized Laplacian
$\D_\gn$ and the combinatorial Laplacian $\D$ are related by the
simple identity $\D=\vk_+\D_\gn$. However, in the case of an
arbitrary graph the spectra of these operators, in spite of many
similar properties, may have significant differences.
\end{remark}

The normalized Laplacian $\D_\gn$ has the decomposition into a
constant fiber direct integral \er{raz}. The precise expression  of
the fiber Laplacian $\D_\gn(k)$ is given by \er{fNL1}, \er{ftro}. The
eigenvalues of the fiber Laplacian $\D_\gn(k)$ will be denoted by
$\m_j(k)$, $j\in\N_\n$. The spectral bands $\s_j(\D_\gn)$,
$j\in\N_{\n}$, for the normalized Laplace operator $\D_\gn$ have the
form $\s_j(\D_\gn)=[\m_j^{-},\m_j^{+}]=\m_j(\T^d)$.

Now we estimate the total bandwidth for the normalized Laplacian $\D_\gn$.

\begin{theorem}\lb{TVENL}
Let $\D_\gn$ be the normalized Laplacian defined by \er{DNLA} on a
periodic graph $\cG$. Then its total bandwidth $\gS(\D_\gn)$
satisfies
\[\lb{CEN111}
\frac{2d_*}{\vk_+^\nu}\le \gS(\D_\gn)\leq\frac{4\cI}{\vk_-}\,,
\]
where the invariant $\cI$ is defined by \er{dIm}, $d_*$ is given in
\er{CEC11},  and $\vk_\pm$ have the form \er{vl+-}.
\end{theorem}

\begin{remarks}
1) The upper estimate in \er{CEN111} is a direct consequence of the estimate \er{esNL} proved in \cite{KS19} and the definition \er{dIm} of the invariant $\cI$. The proof of the lower estimate in \er{CEN111} is similar to the proof of the lower estimate in \er{CEC1}.

2) The lower estimates in \er{CEC1}, \er{CED1} and \er{CEN111}  are
expressed in terms of only the dimension $d$ of the periodic graph,
the number $\n$ of fundamental graph vertices, vertex degrees, and
the potential $V$.
\end{remarks}

\subsection{Examples}
In this subsection we apply the obtained estimates of the total bandwidth to the Laplacians on some simple periodic graphs. In addition we compute the total bandwidth.

\setlength{\unitlength}{1.0mm}
\begin{figure}[h]
\centering
\unitlength 1mm 
\linethickness{0.4pt}
\ifx\plotpoint\undefined\newsavebox{\plotpoint}\fi 

\begin{picture}(160,55)(0,0)
\bezier{25}(21,15)(26,25)(31,35)
\bezier{25}(22,15)(27,25)(32,35)
\bezier{25}(23,15)(28,25)(33,35)
\bezier{25}(24,15)(29,25)(34,35)
\bezier{25}(25,15)(30,25)(35,35)
\bezier{25}(26,15)(31,25)(36,35)
\bezier{25}(27,15)(32,25)(37,35)
\bezier{25}(28,15)(33,25)(38,35)
\bezier{25}(29,15)(34,25)(39,35)
\bezier{25}(30,15)(35,25)(40,35)
\bezier{25}(31,15)(36,25)(41,35)
\bezier{25}(32,15)(37,25)(42,35)
\bezier{25}(33,15)(38,25)(43,35)
\bezier{25}(34,15)(39,25)(44,35)
\bezier{25}(35,15)(40,25)(45,35)
\bezier{25}(36,15)(41,25)(46,35)
\bezier{25}(37,15)(42,25)(47,35)
\bezier{25}(38,15)(43,25)(48,35)
\bezier{25}(39,15)(44,25)(49,35)

\put(8.0,51){$\bK$}
\put(-3,15){\line(1,0){66.0}}
\put(7,35){\line(1,0){66.0}}
\put(17,55){\line(1,0){66.0}}
\put(-1.5,12){\line(1,2){23.0}}
\put(18.5,12){\line(1,2){23.0}}
\put(38.5,12){\line(1,2){23.0}}
\put(58.5,12){\line(1,2){23.0}}

\put(11.5,12){\line(-1,2){8.0}}
\put(31.5,12){\line(-1,2){18.0}}
\put(51.5,12){\line(-1,2){23.0}}
\put(66.5,22){\line(-1,2){18.0}}
\put(76.5,42){\line(-1,2){8.0}}

\put(20,15){\vector(1,0){20.0}}
\put(20,15){\vector(1,2){10.0}}

\put(0,15){\circle{1.5}}
\put(10,15){\circle{1.5}}
\put(20,15){\circle*{1.5}}
\put(30,15){\circle*{1.5}}
\put(40,15){\circle{1.5}}
\put(50,15){\circle{1.5}}
\put(60,15){\circle{1.5}}

\put(5,25){\circle{1.5}}
\put(25,25){\circle*{1.5}}
\put(45,25){\circle{1.5}}
\put(65,25){\circle{1.5}}

\put(10,35){\circle{1.5}}
\put(20,35){\circle{1.5}}
\put(30,35){\circle{1.5}}
\put(40,35){\circle{1.5}}
\put(50,35){\circle{1.5}}
\put(60,35){\circle{1.5}}
\put(70,35){\circle{1.5}}

\put(15,45){\circle{1.5}}
\put(35,45){\circle{1.5}}
\put(55,45){\circle{1.5}}
\put(75,45){\circle{1.5}}

\put(20,55){\circle{1.5}}
\put(30,55){\circle{1.5}}
\put(40,55){\circle{1.5}}
\put(50,55){\circle{1.5}}
\put(60,55){\circle{1.5}}
\put(70,55){\circle{1.5}}
\put(80,55){\circle{1.5}}

\put(32.0,12.5){$\scriptstyle a_1$}
\put(24.0,31.0){$\scriptstyle a_2$}
\put(16.1,12.2){$\scriptstyle O=x_1$}
\put(22.1,36.0){$\scriptstyle x_1+a_2$}
\put(39.8,36.0){$\scriptstyle x_3+a_2$}
\put(40.0,12.5){$\scriptstyle x_1+a_1$}
\put(20.8,25.0){$\scriptstyle x_2$}
\put(46.1,24.5){$\scriptstyle x_2+a_1$}
\put(49.0,32.0){$\scriptstyle x_1+a_1+a_2$}
\put(27.5,12.5){$\scriptstyle x_3$}
\put(33,24){$\Omega$}
\put(-10,15.0){\emph{a})}

\put(125.0,50){$\bK_*$}
\put(98,27){$\Omega$}
\multiput(94,41)(5,0){5}{\line(1,0){3}}
\bezier{30}(107.5,18)(108.5,20)(109.5,22)
\bezier{30}(111.5,26)(112.5,28)(113.5,30)
\bezier{30}(115.5,34)(116.5,36)(117.5,38)

\put(96.5,16.3){$\be_6$}
\put(85.0,16.3){$\be_3$}
\put(90.5,21){$\be_2$}
\put(78.5,21){$\be_1$}
\put(90.0,32){$\be_4$}
\put(105,33.5){$\be_5$}

\put(80,15){\vector(1,0){26.0}}
\put(80,15){\vector(1,2){13.0}}

\bezier{30}(106,15)(112.5,28)(119,41)

\bezier{30}(81,15)(87.5,28)(94,41)
\bezier{30}(82,15)(88.5,28)(95,41)
\bezier{30}(83,15)(89.5,28)(96,41)
\bezier{30}(84,15)(90.5,28)(97,41)
\bezier{30}(85,15)(91.5,28)(98,41)
\bezier{30}(86,15)(92.5,28)(99,41)
\bezier{30}(87,15)(93.5,28)(100,41)
\bezier{30}(88,15)(94.5,28)(101,41)
\bezier{30}(89,15)(95.5,28)(102,41)
\bezier{30}(90,15)(96.5,28)(103,41)
\bezier{30}(91,15)(97.5,28)(104,41)
\bezier{30}(92,15)(98.5,28)(105,41)
\bezier{30}(93,15)(99.5,28)(106,41)
\bezier{30}(94,15)(100.5,28)(107,41)
\bezier{30}(95,15)(101.5,28)(108,41)
\bezier{30}(96,15)(102.5,28)(109,41)
\bezier{30}(97,15)(103.5,28)(110,41)
\bezier{30}(98,15)(104.5,28)(111,41)
\bezier{30}(99,15)(105.5,28)(112,41)
\bezier{30}(100,15)(106.5,28)(113,41)
\bezier{30}(101,15)(107.5,28)(114,41)
\bezier{30}(102,15)(108.5,28)(115,41)
\bezier{30}(103,15)(109.5,28)(116,41)
\bezier{30}(104,15)(110.5,28)(117,41)
\bezier{30}(105,15)(111.5,28)(118,41)

\put(80,15){\circle*{1.5}}
\put(93,15){\circle*{1.5}}
\put(106,15){\circle{1.5}}

\put(93,15){\line(-1,2){6.5}}
\put(112.5,28){\line(-1,2){6.5}}

\put(86.5,28){\circle*{1.5}}
\put(112.5,28){\circle{1.5}}

\put(93,41){\circle{1.5}}
\put(106,41){\circle{1.5}}
\put(119,41){\circle{1.5}}

\put(96.7,12.0){$a_1$}
\put(84.5,34.5){$a_2$}
\put(78.0,11.5){$x_1$}
\put(91.0,42.5){$x_1$}
\put(102.0,42.7){$x_3$}
\put(104.0,11.5){$x_1$}
\put(81.0,28.0){$x_2$}
\put(113.8,28.0){$x_2$}
\put(116.0,42.5){$x_1$}
\put(90.0,11.5){$x_3$}
\put(70,15.0){\emph{b})}

\put(125,15){\line(1,0){26.0}}
\put(125,15){\line(1,2){13.0}}
\put(151,15){\line(-1,2){13.0}}
\put(125,15){\circle*{1.5}}
\put(120.0,13.5){$x_1$}
\put(151,15){\circle*{1.5}}
\put(152.5,13.5){$x_3$}
\put(138,41){\circle*{1.5}}
\put(138.5,42){$x_2$}
\bezier{200}(125,15)(138,3)(151,15)
\bezier{200}(125,15)(120,33)(138,41)
\bezier{200}(151,15)(156,33)(138,41)

\put(136.5,16){$\be_3$}
\put(136.5,5.5){$\be_6$}
\put(132.5,26){$\be_1$}
\put(140,26){$\be_2$}
\put(120.5,30){$\be_4$}
\put(151.5,30){$\be_5$}

\put(131.0,27.0){\vector(1,2){1.0}}
\put(144.2,28.5){\vector(1,-2){1.0}}
\put(138,15.0){\vector(-1,0){1.0}}
\put(150.8,29.5){\vector(-1,2){1.0}}
\put(126.2,31.5){\vector(-1,-2){1.0}}
\put(138,9.0){\vector(1,0){1.0}}
\end{picture}
\vspace{-0.8cm} \caption{\footnotesize   \emph{a}) The Kagome lattice $\bK$; \quad
\emph{b})  the fundamental graph $\bK_*$.} \lb{FEx2}
\end{figure}
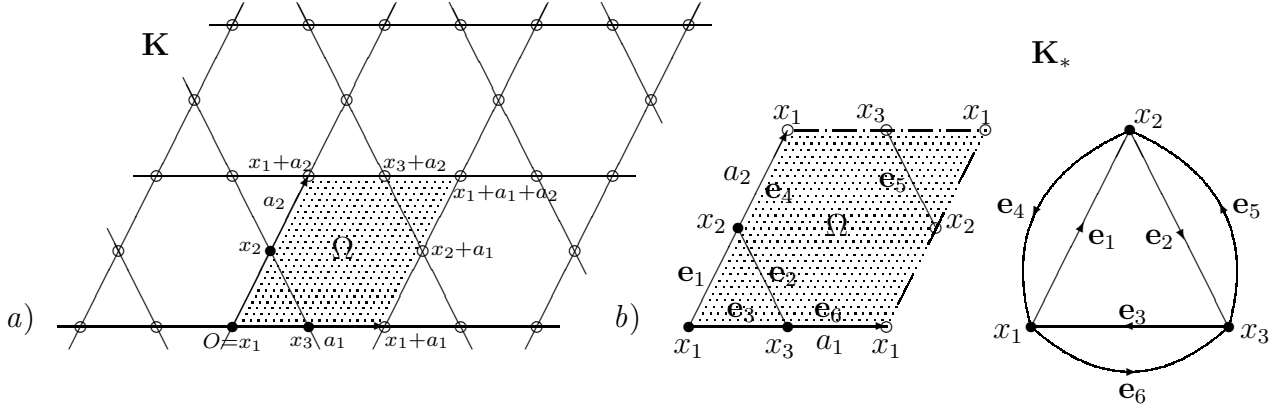

\begin{example}\lb{EKL}
Let $\D$ be the combinatorial Laplacian defined by \er{DLO} on the
Kagome lattice  $\bK$ (see Fig.~\ref{FEx2}a). Then the total
bandwidth $\gS(\D)=6$, but the estimate \er{CED1} is given by
\[\lb{esKL}
\textstyle\frac14\leq \gS(\D)\leq 12.
\]
\end{example}

\begin{remarks}
1) The spectrum of the combinatorial Laplacian $\D$ on the Kagome
lattice is given by $\s(\D)=\s_{ac}(\D)\cup\s_{fb}(\D)$,  where
$\s_{fb}(\D)=\s_3=\{6\}$ and $\s_{ac}(\D)=\s_1\cup\s_2, \
\s_1=[0,3],\  \s_2=[3,6]. $ Then we obtain that the total bandwidth
$\gS(\D)=\sum_{j=1}^3\big|\s_j(\D)\big|=6$.

2)  More accurate lower estimates for the Kagome lattice will be
given in Example \ref{EKL1}.
\end{remarks}

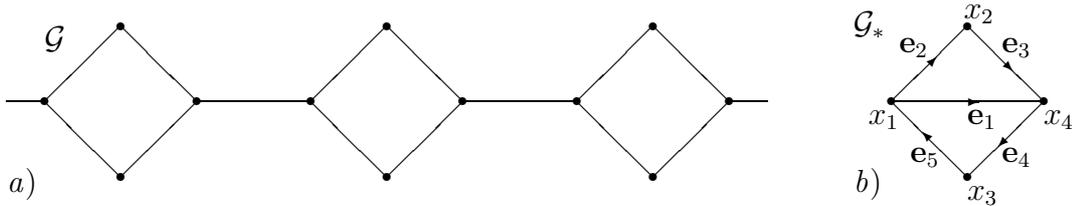
\begin{figure}[h]
\centering
\unitlength 1mm 
\linethickness{0.4pt}
\ifx\plotpoint\undefined\newsavebox{\plotpoint}\fi 
\begin{picture}(100,30)(0,0)

\put(0,15){\line(1,0){5.0}}

\put(15,5){\circle*{1}}
\put(5,15){\circle*{1}}
\put(25,15){\circle*{1}}
\put(15,25){\circle*{1}}
\put(15,5){\line(-1,1){10.0}}
\put(15,5){\line(1,1){10.0}}
\put(15,25){\line(-1,-1){10.0}}
\put(15,25){\line(1,-1){10.0}}

\put(25,15){\line(1,0){15.0}}

\put(50,5){\circle*{1}}
\put(40,15){\circle*{1}}
\put(60,15){\circle*{1}}
\put(50,25){\circle*{1}}
\put(50,5){\line(-1,1){10.0}}
\put(50,5){\line(1,1){10.0}}
\put(50,25){\line(-1,-1){10.0}}
\put(50,25){\line(1,-1){10.0}}

\put(60,15){\line(1,0){15.0}}

\put(85,5){\circle*{1}}
\put(75,15){\circle*{1}}
\put(95,15){\circle*{1}}
\put(85,25){\circle*{1}}
\put(85,5){\line(-1,1){10.0}}
\put(85,5){\line(1,1){10.0}}
\put(85,25){\line(-1,-1){10.0}}
\put(85,25){\line(1,-1){10.0}}

\put(95,15){\line(1,0){5.0}}

\put(0,3){\emph{a})}
\put(5,22){$\cG$}

\end{picture}\hspace{10mm}
\begin{picture}(30,30)(0,0)
\put(15,5){\circle*{1}}
\put(5,15){\circle*{1}}
\put(25,15){\circle*{1}}
\put(15,25){\circle*{1}}
\put(10,10){\vector(-1,1){1.0}}
\put(20,10){\vector(-1,-1){1.0}}
\put(15,5){\line(-1,1){10.0}}
\put(15,5){\line(1,1){10.0}}
\put(15,25){\line(-1,-1){10.0}}
\put(15,25){\line(1,-1){10.0}}
\put(10,20){\vector(1,1){1.0}}
\put(20,20){\vector(1,-1){1.0}}
\put(5,15){\line(1,0){20.0}}
\put(15.5,14.9){\vector(1,0){1.0}}
\put(0,3){\emph{b})}
\put(0,24){$\cG_*$}
\put(6,21.5){$\be_2$}
\put(19.5,21.5){$\be_3$}
\put(7.5,7){$\be_5$}
\put(19.5,7.0){$\be_4$}
\put(2,12.0){$x_1$}
\put(25,12.0){$x_4$}
\put(15,2.0){$x_3$}
\put(14.5,26.0){$x_2$}
\put(15,12.0){$\be_1$}
\end{picture}
\caption{\footnotesize  \emph{a}) A $\Z$-periodic graph $\cG$; \quad \emph{b}) the
fundamental graph  $\cG_*$.} \label{FEx1}
\end{figure}

\begin{example}\lb{E1D}
Let $\D_\gn$ be the normalized Laplacian defined by \er{DNLA} on the
periodic graph $\cG$ shown in Fig.~\ref{FEx1}a. Then the total
bandwidth $\gS(\D_\gn)=\frac43$,  but the estimate \er{CEN111} is
given by
\[\lb{es1D}
\textstyle \frac4{81}\leq\gS(\D_\gn)\leq2 .
\]
\end{example}

\begin{remarks}
1) The spectrum of the normalized Laplacian $\D_\gn$ on the graph
$\cG$ shown in Fig.~\ref{FEx1}\emph{a} is given by
$\s(\D_\gn)=\s_{ac}(\D_\gn)\cup\s_{fb}(\D_\gn)$, where $
\s_{ac}(\D_\gn)=\s_1\cup\s_2\cup\s_3$, and $\s_{fb}(\D_\gn)=\s_4=\{1\}$,
$\s_1=\big[0,\frac13\,\big]$, $\s_2=\big[\frac23\,,\frac43\,\big]$,  $\s_3=\big[\frac53\,,2\big]$.
Then we obtain that the total bandwidth $\gS(\D_\gn)=\frac43$.

2)  More accurate lower estimates for this graph $\cG$ will be given
in Example \ref{E1D1}.
\end{remarks}

\section{Trace formulas for Schr\"odinger operators} \lb{Sec3}
\setcounter{equation}{0}

\subsection{Trace formulas for adjacency operators} We consider the
\emph{adjacency} operator $A$ acting on the Hilbert space $\ell^2(\cV)$ and given by
\[
\lb{AdO}
(A f)_x=\sum_{(x,y)\in\cA}f_y, \qqq f\in\ell^2(\cV), \qqq x\in\cV.
\]
The fiber \emph{adjacency} operator $A(k)$, $k\in\T^d$, has the form
\[
\label{fado}
\big(A(k)f\big)_x=\sum_{\be=(x,y)\in\cA_*} e^{i\lan\t(\be),\,k\ran}f_y,
 \qqq f\in\ell^2(\cV_*),\qqq x\in \cV_*,
\]
where $\t(\be)$ is the edge index, defined by \er{in}, \er{dco}. We
recall Theorem 2.2 from \cite{KS21}.

\begin{theorem}\lb{TFao} Let $A(k)$, $k\in\T^d$, be the fiber adjacency operator defined by \er{fado} on the
fundamental  graph $\cG_*$. Let $\cC_n$ be the set of all cycles of
length $n$ in $\cG_*$, and $\t(\bc)$ be the index of $\bc$ defined
by \er{cyin}. Then for each $n\in\N$ the trace of $A^n(k)$ satisfies
$$
\Tr A^n(k)=\sum_{\bc\in\cC_n}\cos\lan \t(\bc),k\ran,
$$
$$
\frac1{(2\pi)^d}\int_{\T^d}\Tr A^n(k)dk=\cN_n^0\ge 0,
$$
where $\cN_n^0$ is the number of all cycles of length $n$ with zero index in $\cG_*$.
\end{theorem}

\subsection{Trace formulas for Schr\"odinger operators}
In order to formulate trace formulas for the Schr\"odinger operator
$H=-\D+V$  defined by \er{Sh} -- \er{DLO}, we need some modification
of the fundamental graph $\cG_*$. We add a loop $\be_x$ with index
$\t(\be_x)=0$ at each vertex $x$ of the fundamental graph
$\cG_*=(\cV_*,\cA_*)$ and consider the modified fundamental graph
$\wt\cG_*=(\cV_*,\wt\cA_*)$, where
$$
\wt\cA_*=\cA_*\cup\{\be_x\}_{x\in\cV_*}.
$$
We denote by $\wt\cC$ the set of all cycles in $\wt\cG_*$.
For each cycle $\bc\in\wt\cC$ we define the \emph{weight}
\[\lb{Wcy}
\o(\bc)=\o(\be_1)\ldots \o(\be_n), \qq \textrm{where}\qq \bc=(\be_1,\ldots,\be_n)\in\wt\cC,
\]
and $\o(\be)$ is given by
\[\lb{webe}
\o(\be)=\left\{
\begin{array}{cl}
1,  & \qq \textrm{if} \qq  \be\in\cA_* \\[2pt]
v_x, & \qq \textrm{if} \qq \be=\be_x
\end{array}\right.,\qqq v_x=V_x-\vk_x,
\]
where $\vk_x$ is the degree of the vertex $x$.

\begin{remark}
Note that
$$
\o(\bc)=1\qqq \textrm{for each cycle}\qqq \bc\in\cC.
$$
\end{remark}

We recall Theorems 2.4 and 2.5 from \cite{KS21}.

\begin{theorem}\lb{TPG}
Let $H(k)$, $k\in\T^d$, be the fiber Schr\"odinger operator defined
by  \er{Hk'} -- \er{fLao} on the fundamental graph $\cG_*$. Let
$\wt\cC_n$ be the set of all cycles of length $n$ in the modified
fundamental graph $\wt\cG_*$. Then for each $n\in\N$ the trace of $H^n(k)$ satisfies
\[\lb{TrnH}
\Tr H^n(k)=\cT_n(k), \qqq
\cT_n(k)=\sum_{\bc\in\wt\cC_n}\o(\bc)\cos\lan\t(\bc),k\ran,
\]
\[\lb{ITrH}
\frac1{(2\pi)^d}\int_{\T^d}\Tr H^n(k)dk=\cT_{n,0},\qqq
\cT_{n,0}=\sum_{\bc\in\wt\cC_n^0}\o(\bc),
\]
where $\t(\bc)$ is the index of $\bc$ defined by \er{cyin};
$\o(\bc)$ is given by \er{Wcy}, and $\wt\cC_n^0$ is the set of all cycles of length $n$ with zero index in $\wt\cG_*$.
\end{theorem}

\begin{remarks}
1) The formulas \er{TrnH}, \er{ITrH} are \emph{trace formulas}, where the
traces of  the fiber operators are expressed in terms of some
geometric parameters of the graph (vertex degrees, cycle indices and
lengths) and the potential $V$.

2) The index of the cycle $\bc$ in the fundamental graph $\cG_*$ is equal to zero if and only if $\bc$ corresponds to a cycle in the periodic graph $\cG$, see Remark 1 after Fig. \ref{FGHei}.

3) We sometimes write $\cT_n(k,V),\o(\bc,V),\ldots\,$ instead  of
$\cT_n(k),\o(\bc),\ldots\,$, when several potentials $V$ are dealt
with.

4) The trace formulas for the unperturbed operator $H_o=-\D$,  where
$\D$ is the combinatorial Laplacian defined by \er{DLO} are given by
the identities \er{TrnH} and \er{ITrH}, where $\o(\bc)=\o(\bc,0)$.
\end{remarks}

\section{Estimates for the combinatorial Schr\"odinger operators}\lb{Sec4}
\setcounter{equation}{0}

In this section we estimate the total bandwidth for  the
combinatorial Schr\"odinger operators with periodic potentials.

\subsection{Estimates for adjacency operators} Let $\cC_n$ be the set of all cycles of length $n$ in the fundamental graph $\cG_*$. We define the following subsets of $\cC_n$:

$\bu$  $\cC_n^+$ is the set of all cycles of length $n$ with
non-zero indices,  and  $\cN_n^+$ is their number:
\[\lb{cN+}
\cC_n^+=\{\bc\in\cC_n : \t(\bc)\neq0\}, \qqq \cN_n^+=\#\cC_n^+<\iy.
\]

$\bu$  $\cC_n^{odd}$ is the set of all cycles of length $n$ with odd
sum  of index components, and  $\cN_n^{odd}$ is their number:
\[\lb{cNo}
\cC_n^{odd}=\{\bc\in\cC_n : \lan\t(\bc),\1\ran\textrm{ is odd}\}, \qquad \cN_n^{odd}=\#\cC_n^{odd}<\iy,
\]
where $\t(\bc)$ is the index of $\bc$ defined by \er{cyin}, and
$\1=(1,\ldots,1)\in\R^d$.

\begin{remark}
Note that $\cN_{n}^+\geq\cN_{n}^{odd}$ for each $n\in\N$, since $\cC_{n}^{odd}\subseteq\cC_{n}^+$.
\end{remark}

We consider the adjacency operator $A$.

\begin{theorem}\lb{TeAO}
Let $A$ be the adjacency operator defined by \er{AdO} on a periodic graph  $\cG$, and let $n\in\N$. Then the total bandwidth $\gS(A^n)=\sum_{j=1}^{\n}|\s_j(A^n)|$ for the operator $A^n$ satisfies
$$
\gS(A^n)\geq \max\{\cN_n^+,2\cN_n^{odd}\},
$$
where $\cN_n^+$ and $\cN_n^{odd}$ are defined in \er{cN+} and \er{cNo}.
\end{theorem}

We omit the proof, since it is similar to the proof of Theorem \ref{TEsS}.

\begin{remark}
The numbers $\cN_n^+$ and $\cN_n^{odd}$ defined by \er{cN+} and \er{cNo} are
the numbers of all cycles from the corresponding sets
including cycles with \emph{back-tracking parts}, i.e., including
cycles $(\be_1,\ldots,\be_n)$ for which $\be_{s+1}=\ul\be_s$ for
some $s\in\N_n$ ($\be_{n+1}$ is understood as $\be_1$).
\end{remark}

Now we estimate the total bandwidth
$\gS(A)=\sum_{j=1}^\n\big|\s_j(A)\big|$ for the adjacency operator
$A$.

\begin{corollary}\lb{CAEs}
Let $A$ be the adjacency operator defined by \er{AdO} on a periodic
graph $\cG$. Then the total bandwidth $\gS(A)$ satisfies
$$
\frac1{n\vk_+^{n-1}}\max\{\cN_n^+,2\cN_n^{odd}\}\leq  \gS(A)  \leq
4\cI, \qqq \forall\,n\in\N,
$$
where the invariant $\cI$ has the form \er{dIm}; $\cN_n^+$  and
$\cN_n^{odd}$ are defined in \er{cN+} and \er{cNo}, and $\vk_+$ is given in
\er{bf}. In particular, the adjacency operator $A$ has at least
one non-degenerate band.
\end{corollary}

We omit the proof, since it is similar to the proof of Corollary \ref{CHEs}.

\subsection{Estimates of bandwidths for Schr\"odinger operators}
\lb{ssSO} We discuss estimates of the total bandwidth for  the
Schr\"odinger operators $H=-\D+V$ on periodic graphs. Without loss
of generality (adding a constant to the periodic potential $V$) we
may assume that
\[
\lb{Qmi0} v_-=\min_{x\in\cV_*}(V_x-\vk_x)=0,\qqq
v_+=\max_{x\in\cV_*}(V_x-\vk_x)\ge 0,
\]
where $\vk_x$ is the degree of the vertex $x$. From $\s(A)\subseteq
[-\vk_+,\vk_+]$ and under the condition \er{Qmi0} we deduce that
$$\s(H)\subseteq[-\vk_+,\vk_++v_+].
$$

\begin{theorem}\lb{TEsS}
Let $H=-\D+V$ be the Schr\"odinger operator defined by \er{Sh} --
\er{DLO} with a periodic potential $V$ satisfying \er{Qmi0} on a
periodic graph $\cG$ with the fundamental graph
$\cG_*=(\cV_*,\cA_*)$, and let $n\in\N$. Then the total bandwidth
$\gS(H^n)$ for the operator $H^n$ satisfies
\[\lb{EsS1}
\gS(H^n)\geq \max\{B_{n,1},B_{n,2}\},
\]
where
\[\lb{coCH}
B_{n,1}=\sum_{\bc\in\wt\cC_{n}\atop\t(\bc)\neq0}\o(\bc)\geq\cN_n^+,\qqq B_{n,2}=
2\hspace{-5mm}\sum_{\bc\in\wt\cC_n\atop \lan\t(\bc),\1\ran\textrm{ is odd}}
\hspace{-5mm}\o(\bc)\geq2\cN_n^{odd}.
\]
Here $\wt\cC_n$ is the set of all cycles of length $n$ in the
modified  fundamental graph $\wt\cG_*$; $\t(\bc)$ is the index of
$\bc$ defined by \er{cyin}; $\o(\bc)$ is given by \er{Wcy};
$\1=(1,\ldots,1)\in\R^d$, and $\cN_n^+$ and $\cN_n^{odd}$ are defined by \er{cN+} and \er{cNo}.
\end{theorem}

\no{\bf Proof}. The definition of the spectral bands of $H^n$ and
the identity \er{TrnH} imply
\begin{multline*}
\gS(H^n)= \sum_{j=1}^\n\big|\s_j(H^n)\big|=
\sum_{j=1}^\n\Big(\max_{k\in\T^d}\l_j^n(k)-
\min_{k\in\T^d}\l_j^n(k)\Big)\geq
\sum_{j=1}^\n\l_j^n(0)-\sum_{j=1}^\n\l_j^n(k_\bu)
\\=\Tr H^n(0)-\Tr H^n(k_\bu)
=\cT_n(0)-\cT_n(k_\bu)\\=\sum_{\bc\in\wt\cC_n}\o(\bc)
\big(1-\cos\lan\t(\bc),k_\bu\ran\big)=2\sum_{\bc\in\wt\cC_n}\o(\bc)
\sin^2\frac{\lan\t(\bc),k_\bu\ran}2
\end{multline*}
for any $k_\bu\in\T^d$. Thus,
\[\lb{kbu}
\gS(H^n)\geq \cT_n(0)-\cT_n(k_\bu)= 2\sum_{\bc\in\wt\cC_n}\o(\bc)
\sin^2\frac{\lan\t(\bc),k_\bu\ran}2\,,\qqq \forall\,k_\bu\in\T^d.
\]

Due to \er{TrnH} and \er{ITrH}, we have
$$
\cT_{n,0}=\frac1{(2\pi)^d}\int_{\T^d}\Tr H^n(k)dk
=\Tr H^n(k_0)=\cT_n(k_0)\qq \textrm{for some}\qq k_0\in\T^d.
$$

If $k_\bu=k_0$, then, using \er{TrnH} and \er{ITrH}, the estimate \er{kbu} has the form
\[\lb{esB1}
\gS(H^n)\geq \cT_n(0)-\cT_n(k_0)=\cT_n(0)-\cT_{n,0}
=\sum_{\bc\in\wt\cC_n\atop\t(\bc)\neq0}\o(\bc).
\]
Due to \er{Wcy}, \er{webe} and \er{Qmi0}, $\o(\bc)\geq0$ for all
$\bc\in\wt\cC$  and $\o(\bc)=1$ for each cycle $\bc\in\cC$. Then we
obtain the following estimate
\[\lb{esBn1}
B_{n,1}=\sum_{\bc\in\wt\cC_n\atop\t(\bc)\neq0}\o(\bc)\geq
\sum_{\bc\in\cC_n\atop\t(\bc)\neq0}\o(\bc)=\cN_n^+.
\]

Similarly, if $k_\bu=\pi\1$, we obtain
\[\label{esB2}
\gS(H^n)\geq 2\sum_{\bc\in\wt\cC_n}\o(\bc)
\sin^2\frac{\lan\t(\bc),\pi\1\ran}2=2\hspace{-5mm}
\sum_{\bc\in\wt\cC_n\atop \lan\t(\bc),\1\ran\textrm{ is
odd}}\hspace{-5mm}\o(\bc)\geq2\hspace{-5mm} \sum_{\bc\in\cC_n\atop
\lan\t(\bc),\1\ran\textrm{ is
odd}}\hspace{-5mm}\o(\bc)=2\cN_n^{odd}.
\]
Here we have also used that $\t(\bc)\in\Z^d$ for all cycles
$\bc\in\wt\cC$.  Combining \er{esB1} -- \er{esB2}, we get \er{EsS1},
\er{coCH}. \qq \BBox

\medskip

Now we estimate the total bandwidth for the Schr\"odinger operator $H$.

\begin{corollary}\lb{CHEs}
Let $H=-\D+V$ be the Schr\"odinger operator defined by \er{Sh} --
\er{DLO}  with a periodic potential $V$ satisfying \er{Qmi0} on a
periodic graph $\cG$. Then its total bandwidth $\gS(H)$ satisfies
\[\lb{CEH1}
\frac{\max\{B_{n,1},B_{n,2}\}}{n\,(v_++\vk_+)^{n-1}}\leq\gS(H)\leq4\cI, \qqq \forall\,n\in\N,
\]
where the invariant $\cI$ is defined by \er{dIm}; $B_{n,s}$, $s=1,2$, have the form \er{coCH}, and $v_+$ and $\vk_+$ are given by \er{vl+-}.
\end{corollary}

\no {\bf Proof}. Recall that under the condition \er{Qmi0} we have
$\s(H)\subseteq [-\vk_+,v_++\vk_+]$, where $\vk_+=\max_{x\in\cV_*}\vk_x$.  We have a
simple estimate for the spectral bands $\s_j(H^n)$ and $\s_j(H)$,
$j=1,\ldots,\n$:
\begin{multline*}
|\s_j(H^n)|=\max_{k\in\T^d}\l_{j}^n(k)-\min_{k\in\T^d}\l_{j}^n(k)=
\l_j^n(k^+)-\l_j^n(k^-)\\\le n(v_++\vk_+)^{n-1}
\big|\l_j(k^+)-\l_j(k^-)\big|\le n(v_++\vk_+)^{n-1}|\s_j(H)|,
\end{multline*}
for some $k^\pm\in\T^d$. Then using \er{EsS1},  we obtain the lower estimate in \er{CEH1}. The upper estimate in \er{CEH1} was proved in \cite{KS20}. \qq \BBox

\medskip

In order to prove Theorem \ref{TVECS} we need Lemma \ref{LFtZ}.
Recall that $\cC$ is the set of all cycles of the fundamental graph
$\cG_*=(\cV_*,\cA_*)$, and for any cycle $\bc\in\cC$ the cycle index
$\t(\bc)\in\Z^d$ is defined by \er{cyin}. Repeating a cycle $\bc$
$m$ times, we obtain \emph{$m$-multiple $\bc^m$} of $\bc$. If $\bc$
is not a $m$-multiple of a cycle with $m\geq2$, $\bc$ is called
\emph{prime}. A cycle $\bc=(\be_1,\ldots,\be_n)$ has
\emph{backtracking}  if $\be_{s+1}=\ul\be_s$, for some $s\in\N_n$
($\be_{n+1}$ is understood as $\be_1$). A cycle with no backtracking
is  called \emph{proper}. We show that in $\cG_*$ there exist $d$
prime proper cycles of length $\le\n=\#\cV_*$ with indices forming
an orthonormal basis of $\Z^d$.

A graph is called \emph{bipartite} if its vertex set is divided into
two  disjoint sets (called \emph{parts} of the graph) such that each
edge connects vertices from distinct parts.

\begin{lemma}\lb{LFtZ}
i) Let $\cC$ be the set of all cycles of the fundamental graph
$\cG_*=(\cV_*,\cA_*)$. Then the image of the function $\t:\cC\ra\R^d$
defined by \er{cyin} is the lattice $\Z^d$, i.e.,
\[\lb{prao}
\t(\cC\,)=\Z^d.
\]

ii) For some basis of the lattice $\G$ there exist prime cycles
$\bc_1,\ldots,\bc_d$ of lengths $\le \n=\#\cV_*$, such that the set
$\big\{\t(\bc_s)\big\}_{s\in\N_d}$ forms an orthonormal basis of
$\Z^d$.

iii) There exists $n\leq\n$ such that the number $\cN_n^{odd}$ defined by \er{cNo} satisfies
\[\lb{cNpo}
\cN_{n}^{odd}\geq nd_*,\qqq
d_*=\left\{\begin{array}{cl}
              d, & \textrm{if $d$ is even} \\
              d+1, & \textrm{if $d$ is odd}
            \end{array}\right..
\]
If $\cG$ is bipartite, then
\[\lb{cNpob}
\cN_{n}^{odd}\geq 2nd.
\]
\end{lemma}

\no \textbf{Proof.} \emph{i}) Let $\G\ss\R^d$ be a lattice with  a
basis $\A=\{a_1,\ldots,a_d\}$. We consider a $\G$-periodic graph
$\cG=(\cV,\cA)$. We take some vertex $x\in\cV$ of $\cG$. Let
$\mm=(m_j)_{j=1}^d\in\Z^d$. Then
$a=\sum\limits_{s=1}^dm_sa_s\in\G$. Since $\cG$ is connected, there
exists an oriented path $\bp$ from $x$ to $x+a$ in $\cG$. Then
$\bc=\bp/\G$ is a cycle in the fundamental graph $\cG_*=\cG/\G$ with
index $\t(\bc)=\mm$. Thus, $\mm\in\t(\cC)$. Conversely, let
$\mm=(m_j)_{j=1}^d\in\R^d$ and $\mm=\t(\bc)$ for some cycle
$\bc\in\cC$. Then there exists an oriented path $\bp$ in the
periodic graph $\cG$ such that $\bc=\bp/\G$ and $\bp$ connects
$\G$-equivalent vertices $x$ and $x+a$ in $\cG$, where
$a=\sum\limits_{s=1}^dm_sa_s$. Thus, $a\in\G$ and, consequently,
$\mm\in\Z^d$.

\emph{ii}) Let $\cT=(\cV_*,\cE_\cT)$ be a spanning tree of  the
fundamental graph $\cG_*$. Then for each
$\be\in\cS_\cT=\cE_*\sm\cE_\cT$ there exists a cycle $\bc_\be$ in
$\cG_*$ whose edges are all in $\cT$ except $\be$. Each $\bc_\be$ is
prime and proper and the set of all cycles $\cB=\{\bc_\be\}_{\be\in\cS_\cT}$
forms a basis of the cycle space $\cC$ of the graph $\cG_*$. Note
that the length of each cycle from $\cB$ is not more than $\n$.
Since $\t(\cC)=\Z^d$ (see \er{prao}), we conclude that
$\{\t(\bc)\}_{\bc\in\cB}$ generates the group $\Z^d$. Consequently,
there exist $d$ cycles $\bc_1,\ldots,\bc_d\in\cB$ such that the set
$$
\big\{\t(\bc_s)=\big(\t_1(\bc_s),\ldots,\t_d(\bc_s)\big)\in\Z^d:
s\in\N_d\big\}
$$
forms a basis of $\Z^d$. Then the set
$$
\A'=\big\{a'_s=\t_1(\bc_s)a_1+\ldots+\t_d(\bc_s)a_d: s\in\N_d\big\}
$$
is a basis of the lattice $\G$. Thus, under the change of the basis
of the lattice $\G$ from $\A$ to $\A'$ the indices
$\big\{\t(\bc_s)\big\}_{s\in\N_d}$ of the cycles
$\bc_1,\ldots,\bc_d$ transform to an orthonormal basis of $\Z^d$.

\emph{iii}) We divide the set $\{\bc_1,\ldots,\bc_d\}$ into two
disjoint sets:  the set $\cS_e$ of cycles with even length and the
set $\cS_o$ of cycles with odd length. Denote by $p$ the maximum
number of entries in the sets $\cS_e$ and $\cS_o$:
\[\lb{esp}
p=\max\{\#\cS_e,\#\cS_o\}, \qqq p\geq\frac{d_*}2\,,
\]
where $d_*$ is given in \er{cNpo}. Let $p=\#\cS_e$. The proof for
the case $p=\#\cS_o$ is similar. Denote by $n$ the maximal length
of cycles from the set $\cS_e$. If some cycle $\bc\in\cS_e$ has the
length less than $n$, then we consider the cycle
$\wt\bc=(\bc,\be,\ul\be,\ldots,\be,\ul\be\,)\in\cC$ of length $n$
with a \emph{back-tracking part} $(\be,\ul\be,\ldots,\be,\ul\be\,)$
for some $\be\in\cA_*$. Due to item \emph{ii}), $n\leq\n$, and the
cycle $\bc$ is prime and proper. Then the cycle $\wt\bc$ is also
prime and $\t(\wt\bc\,)=\t(\bc)$, since $\t(\ul\be\,)=-\t(\be)$ for
all $\be\in\cA_*$.

Using that each cyclic permutation of edges of a prime cycle $\bc$
gives a prime cycle of the same length $|\bc|$ and with the same
index $\t(\bc)$, and the reverse $\ul\bc$ of a prime cycle $\bc$ is
also a prime cycle with $|\ul\bc|=|\bc|$ and $\t(\ul\bc)=-\t(\bc)$, we
obtain
\[\lb{cNgr}
\cN_{n}^{odd}\geq2np.
\]
Combining this and \er{esp}, we obtain \er{cNpo}.

If $\cG$ is bipartite, then without loss of generality we may assume
that  the fundamental graph $\cG_*$ is also bipartite. Then there
are no cycles of odd length in $\cG_*$, and $p=\#\cS_e=d$. Combining
this and \er{cNgr}, we obtain \er{cNpob}. \qq $\BBox$

\medskip

Now we prove Theorem \ref{TVECS} and Corollaries \ref{CVECS}, \ref{ESO}.

\medskip

\no \textbf{Proof of Theorem \ref{TVECS}.} By Lemma
\ref{LFtZ}.\emph{iii}), there exists $n\leq\n$ such that  the number
$\cN_n^{odd}$ defined by \er{cNo} satisfies $\cN_{n}^{odd}\geq nd_*$,
where $d_*$ is given in \er{CEC11}. Then, using the second
inequality in  \er{coCH}, \er{CEH1} and $v_*=v_++\vk_+\geq2$ we
have
$$
\gS(H)=\sum_{j=1}^\n\big|\s_j(H)\big|\geq\frac{B_{n,2}}{n\,v_*^{n-1}}\geq
\frac{2\cN_{n}^{odd}}{n\,v_*^{n-1}}\geq
\frac{2nd_*}{n\,v_*^{n-1}}\geq\frac{2d_*}{v_*^{\n-1}}\,.
$$
The upper estimate in \er{CEC1} was proved in \cite{KS20}.
The estimate \er{LMH} follows from \er{CEC1} and the inequalities
$$
\frac1\n\,\gS(H)\leq\max_{j\in\N_\n}|\s_j(H)|\leq|\s(H)|\leq\gS(H)=
\sum_{j=1}^{\n}|\s_j(H)|. \qqq \BBox
$$

\medskip

\no \textbf{Proof of Corollary \ref{CVECS}.} The estimate \er{CED1}
for the operator $H_o$ is obtained from the estimate \er{CEC1} for
the Schr\"odinger operator $H=H_o+V$ as $V=0$. \qq \BBox

\begin{remark}
If $\cG$ is bipartite, then, by \er{cNpob}, $\cN_{n}^{odd}\geq 2nd$
for some $n\leq\n$ and the lower estimates in \er{CEC1} and
\er{CED1} can be improved:
$$
\gS(H)\geq\frac{2\cN_{n}^{odd}}{n\,v_*^{n-1}}
\geq
\frac{4nd}{n\,v_*^{n-1}}\geq\frac{4d}{v_*^{\n-1}}\,,\qqq
\gS(H_o)\geq\frac{4d}{\vk_*^{\n-1}}\,.
$$
\end{remark}

\no \textbf{Proof of Corollary \ref{ESO}.} It is well known (see,
e.g., \cite{vM76}) that the spectrum of  the one-dimensional
Schr\"odinger  operator $H=-\D+V$ on $\Z$ consists of $\n$ bands
$$
\s_n=[\l_n^-,\l_n^+],\qq n\in\N_\n,\qqq \l_1^-<\l_1^+\le
\l_2^-<\ldots<\l_{\n-1}^+\le \l_\n^-<\l_\n^+,
$$
where $\l_\n^+,\l_{\n-1}^-,\l_{\n-2}^+,\ldots$ are the eigenvalues
of the fiber operator  $H(0)$;
$\l_\n^-,\l_{\n-1}^+,\l_{\n-2}^-,\ldots$ are the eigenvalues of
$H(\pi)$ (see Fig. \ref{fig}). These bands are separated by gaps
$(\l_n^{+},\l_{n+1}^{-})$, $n\in\N_{\n-1}$. Some of the gaps may be
degenerate, i.e., $\l_n^{+}=\l_{n+1}^{-}$.

For the one-dimensional lattice  $\Z$, the numbers $d_*$ and $v_*$
defined in \er{CEC11} have the form
$$
v_*=\vk_++\diam(V-\vk)=2+\diam V,\qqq d_*=2.
$$
The fundamental graph $\cG_*$ of the lattice $\Z$ is just the
cycle graph with $\n$ vertices $1,\ldots,\n$, and $\n$ edges
$$
\be_1=(1,2),\qq \be_2=(2,3),\qq \ldots, \qq \be_{\n-1}=(\n-1,\n),\qq \be_\n=(\n,1)
$$
with indices
$$
\t(\be_1)=\ldots=\t(\be_{\n-1})=0, \qqq \t(\be_\n)=1,
$$
and their inverse edges. Thus, $\cI=d=1$, and the estimate \er{CEC1} has the
form \er{1diS}.\qq \BBox

\setlength{\unitlength}{1.0mm}
\begin{figure}[h]
\centering
\unitlength 1.0mm 
\begin{picture}(135,20)
\put(5,10){\line(1,0){120.00}}
\linethickness{3pt} \put(10,10.5){\line(1,0){15.00}}
\put(35,10.5){\line(1,0){10.00}} \put(57,10.5){\line(1,0){17.00}}
\put(85,10.5){\line(1,0){12.00}} \put(106,10.5){\line(1,0){15.00}}
\linethickness{0.4pt} \put(16,12.5){$\s_1$} \put(9,6){$\l_1^-$}
\put(23,6){$\l_1^+$} \put(38,12.5){$\s_2$} \put(34,6){$\l_2^-$}
\put(42,6){$\l_2^+$} \put(49,6){$\ldots$} \put(61,12.5){$\s_{\n-2}$}
\put(56,6){$\l_{\n-2}^-$} \put(71,6){$\l_{\n-2}^+$}
\put(87,12.5){$\s_{\n-1}$} \put(84,6){$\l_{\n-1}^-$}
\put(95,6){$\l_{\n-1}^+$} \put(112,12.5){$\s_\n$}
\put(105,6){$\l_\n^-$} \put(119,6){$\l_\n^+$}
\end{picture}
\vspace{-4mm} \caption{\footnotesize The spectrum of the
Schr\"odinger operator $H=-\D+V$ with a $\n$-periodic potential $V$ on $\Z$.} \label{fig}
\end{figure}
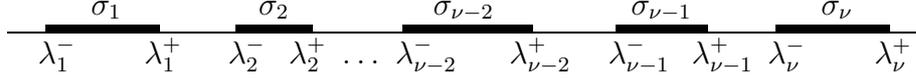

\section{Estimates for normalized Laplacians}\lb{Sec5}
\setcounter{equation}{0}
\subsection{Direct integral decomposition}
We recall Proposition 1.1 from \cite{KS19}.

\begin{theorem}\label{TFD1}
Let the Hilbert space $\mH$ be defined by \er{Hisp}. Then the
normalized  Laplacian $\D_\gn$ on $\ell^2(\cV)$ given by \er{DNLA}
has the following decomposition into a constant fiber direct
integral
\[
\lb{raz}
\begin{aligned}
& U\D_\gn U^{-1}=\int^\oplus_{\T^d}\D_\gn(k)\,{dk\/(2\pi)^d}\,,
\end{aligned}
\]
for some unitary operator $U:\ell^2(\cV)\to\mH$ (the Gelfand
transform).  The fiber Laplacian $\D_\gn(k)$, $k\in\T^d$, on
$\ell^2(\cV_*)$ is given by
\[
\lb{fNL1}
\D_\gn(k)=I-T(k).
\]
Here $I$ is the identity operator on $\ell^2(\cV_*)$, and the fiber
transition  operator $T(k)$ has the form
\[
\label{ftro}
\big(T(k)f\big)_x=\sum_{\be=(x,y)\in\cA_*}
\frac{e^{i\lan\t(\be),\,k\ran}}{\sqrt{\vk_x\vk_y}}\,f_y, \qqq f\in\ell^2(\cV_*),
\qqq x\in \cV_*,
\]
where $\vk_x$ is the degree of the vertex $x$, and $\t(\be)$ is the
index of  the edge $\be\in\cA_*$ defined by \er{in}, \er{dco}.
\end{theorem}

\subsection{Trace formulas for transition operators} Recall that $\cC$ is
the set of all cycles in the fundamental graph
$\cG_*=(\cV_*,\cA_*)$. For each cycle
$\bc=(\be_1,\ldots,\be_n)\in\cC$ we define the \emph{weight}
\[\lb{cywe}
\o_\gn(\bc)=\frac1{\vk_{x_1}\ldots\vk_{x_n}}\,, \qq \textrm{where}\qq
 \be_s=(x_{s},x_{s+1})\in\cA_*, \qq s\in\N_n, \qq x_{n+1}=x_1,
\]
and $\vk_x$ is the degree of the vertex $x$. We recall Theorem 5.3 from \cite{KS21}.

\begin{theorem}\lb{TFNL0} Let $T(k)$, $k\in\T^d$, be the fiber transition
 operator defined by \er{ftro}. Then for each $n\in\N$

i) The trace of $T^n(k)$ has the form
$$
\Tr T^n(k)=T_n(k), \qqq
T_n(k)=\sum_{\bc\in\cC_n}\o_\gn(\bc)\cos\lan\t(\bc),k\ran,
$$
where $\cC_n$ is the set of all cycles of length $n$ in the
fundamental graph $\cG_*$; $\t(\bc)$ is the index of
$\bc$ defined by \er{cyin}, and $\o_\gn(\bc)$ is given by \er{cywe}.

ii) The trace of $T^n(k)$ satisfies
$$
\frac1{(2\pi)^d}\int_{\T^d}\Tr T^n(k)dk=T_{n,0},\qqq
T_{n,0}=\sum_{\bc\in\cC_n^0}\o_\gn(\bc),
$$
where $\cC_n^0$ is the set of all cycles of length $n$ with zero index in $\cG_*$.
\end{theorem}

\subsection{Estimates of bandwidths}
\lb{ssLO} We discuss estimates of the total bandwidth for the normalized Laplacian $\D_\gn$.

\begin{theorem}\lb{PEsN}
Let $T=I-\D_\gn$ be the transition operator, where $\D_\gn$ is  the
normalized Laplacian defined by \er{DNLA} on a periodic graph $\cG$
with the fundamental graph $\cG_*=(\cV_*,\cA_*)$, and let $n\in\N$.
Then the total bandwidth $\gS(T^n)=\sum_{j=1}^\n\big|\s_j(T^n)\big|$, $j\in\N_\n$, $\n=\#\cV_*$, for the operator $T^n$ satisfies
\[\lb{EN1}
\gS(T^n)\geq \max\{B_{n,1},B_{n,2}\},
\]
where
\[\lb{coC}
\begin{aligned}
&\dfrac{\cN_n^+}{\vk_+^n}\leq B_{n,1}=\sum_{\bc\in\cC_{n}\atop\t(\bc)\neq0}\o_\gn(\bc)
\leq\dfrac{\cN_n^+}{\vk_-^n}\,,\\
&\dfrac2{\vk_+^n}\,\cN_n^{odd}\leq B_{n,2}=
2\hspace{-5mm}\sum_{\bc\in\cC_n\atop \lan\t(\bc),\1\ran\textrm{ is odd}}
\hspace{-5mm}\o_\gn(\bc)\leq\dfrac2{\vk_-^n}\,\cN_n^{odd}.
\end{aligned}
\]
Here $\cC_n$ is the set of all cycles of length $n$ in $\cG_*$;
$\t(\bc)$  is the index of $\bc$ defined by \er{cyin}, $\o_\gn(\bc)$
is given by \er{cywe}; $\vk_\pm$ are defined by \er{vl+-}, and the
numbers  $\cN_n^+$ and $\cN_n^{odd}$ are defined in \er{cN+} and
\er{cNo}. For regular graphs of degree $\vk_+$ the inequalities in
\er{coC} become identities.
\end{theorem}

The proof of Theorem \ref{PEsN} is similar to the proof of Theorem \ref{TEsS}.

\medskip

Now we estimate the total bandwidth for the Laplacian $\D_\gn$.

\begin{corollary}\lb{CLEs}
Let $\D_\gn$ be the normalized Laplacian defined by \er{DNLA} on a
periodic  graph $\cG$. Then its total bandwidth $\gS(\D_\gn)$ satisfies
\[\lb{CEN1'}
\frac1n\max\{B_{n,1},B_{n,2}\}\leq\gS(\D_\gn)\leq 2\gb ,
 \qqq \forall\,n\in\N,
\]
where
\[\lb{gb}
\gb=\sum_{x\in \cV_*}\frac{\gb_x}{\vk_x}\leq\frac{2\cI}{\vk_-}\,,
\]
$\gb_x$ is the number of the fundamental graph edges having non-zero
indices and starting at the vertex $x\in\cV_*$; $\vk_x$ is the
degree of $x$; the invariant $\cI$ has the form \er{dIm}; $\vk_-$
is defined in \er{vl+-}, and  $B_{n,s}$, $s=1,2$, are given in \er{coC}.
\end{corollary}

\no {\bf Proof}. Denote by $\x_j(k)$, $j\in\N_\n$, the eigenvalues of the fiber transition operator $T(k)$ defined by \er{ftro}. Then the spectral bands $\s_j(T)$, $j\in\N_{\n}$, for the transition operator $T=I-\D_\gn$ have the form $\s_j(T)=[\x_j^-,\x_j^+]=\x_j(\T^d)$. The spectrum of $T$ satisfies $\s(T)\subseteq[-1,1]$. We have a simple estimate for the spectral bands
$\s_j(T^n)$ and $\s_j(T)$, $j=1,\ldots,\n$:
\[\lb{essb1}
|\s_j(T^n)|=\max_{k\in\T^d}\x_{j}^n(k)-\min_{k\in\T^d}\x_{j}^n(k)=
\x_j^n(k^+)-\x_j^n(k^-)\le n
\big|\x_j(k^+)-\x_j(k^-)\big|\le n|\s_j(T)|,
\]
for some $k^\pm\in\T^d$. Using estimates \er{essb1}, \er{EN1} and the identity
$|\s_j(\D_\gn)|=|\s_{\n-j+1}(T)|$,  $j\in\N_\n$, we obtain the lower estimate in \er{CEN1'}.

The upper estimate in \er{CEN1'} was proved in \cite{KS19}. We prove the inequality in \er{gb}. We have
$$
\gb=\sum_{x\in \cV_*}\frac{\gb_x}{\vk_x}\leq\frac1{\vk_-}\sum_{x\in \cV_*}\gb_x=\frac1{\vk_-}\,\#\cB_*,
$$
where $\vk_-$ is defined in \er{vl+-}, and $\cB_*$ is the set of all bridges of the fundamental graph $\cG_*=(\cV_*,\cA_*)$. This estimate holds for any embedding of the periodic graph $\cG$ into $\R^d$. Choosing the embedding with the minimal number of the fundamental graph bridges, due to \er{dIm}, we obtain
$\gb\leq\frac{2\cI}{\vk_-}$\,. \qq \BBox

\medskip

Now we prove Theorem \ref{TVENL}.

\medskip

\no \textbf{Proof of Theorem \ref{TVENL}.} The proof is similar to
the proof of Theorem \ref{TVECS}. By Lemma \ref{LFtZ}.\emph{iii}),
there exists $n\leq\n$ such that the number  $\cN_n^{odd}$ defined
by \er{cNo} satisfies $\cN_{n}^{odd}\geq nd_*$, where $d_*$ is given
in \er{CEC11}. Then, using \er{CEN1'}, the second inequality in
\er{coC} and $\vk_+\geq2$, we obtain
$$
\gS(\D_\gn)\geq\frac{B_{n,2}}{n}\geq
\frac{2\cN_{n}^{odd}}{n\,\vk_+^n}\geq\frac{2nd_*}{n\,\vk_+^n}
\geq\frac{2d_*}{\vk_+^{\n}}\,.
$$
The upper estimate in \er{CEN111} follows from the upper estimate in \er{CEN1'} and \er{gb}. \qq \BBox

\begin{remark}
If $\cG$ is bipartite, then, by \er{cNpob}, $\cN_{n}^{odd}\geq 2nd$
for  some $n\leq\n$ and the lower estimate in \er{CEN111} can be
improved:
$$
\gS(\D_\gn)\geq\frac{2\cN_{n}^{odd}}{n\,\vk_+^n}\geq
\frac{4nd}{n\,\vk_+^n}\geq\frac{4d}{\vk_+^{\n}}\,.
$$
\end{remark}

\section{Examples} \lb{Sec6}
\setcounter{equation}{0}
First we prove Examples \ref{EKL} and \ref{E1D}
about estimates of bandwidths for the Laplacians on the Kagome
lattice and on the graph shown in Fig.~\ref{FEx1}\emph{a}.

\medskip

\no \textbf{Proof of Example \ref{EKL}.} The Kagome lattice $\bK$ is a
periodic  regular graph of degree $\vk_+=4$. The fundamental graph
$\bK_*$ of $\bK$ consists of three vertices $x_1,x_2,x_3$, six edges
$$
\be_1=(x_1,x_2),\qq \be_2=(x_2,x_3),\qq \be_3=(x_3,x_1),\qq \be_4
=(x_2,x_1),\qq \be_5=(x_3,x_2),\qq \be_6=(x_1,x_3)
$$
with indices
$$
\begin{array}{lll}
\t(\be_1)=(0,0), \qqq & \t(\be_2)=(0,0), \qqq & \t(\be_3)=(0,0),\\[2pt]
\t(\be_4)=(0,1), & \t(\be_5)=(1,-1), & \t(\be_6)=(-1,0),
\end{array}
$$
and their inverse edges, see Fig.~\ref{FEx2}\emph{b}. Thus,
$\cI=3$, $d_*=d=2$, $\vk_*=\vk_+=4$, $\n=3$,  and the estimate \er{CED1} has the form \er{esKL}.\qq \BBox

\medskip

\no \textbf{Proof of Example \ref{E1D}.} We consider the
$\Z$-periodic graph  $\cG$ shown in Fig. \ref{FEx1}\emph{a}. Its
fundamental graph $\cG_*$ consists of four vertices
$x_1,x_2,x_3,x_4$ with degrees $\vk_{x_1}=\vk_{x_4}=3$ and
$\vk_{x_2}=\vk_{x_3}=2$ and five edges
$$
\be_1=(x_1,x_4),\qqq \be_2=(x_1,x_2),\qqq \be_3=(x_2,x_4),\qqq
\be_4=(x_4,x_3), \qqq \be_5=(x_3,x_1)
$$
with indices
$$
\t(\be_1)=1,\qqq \t(\be_2)=\t(\be_3)=\t(\be_4)=\t(\be_5)=0,
$$
and their inverse edges (see Fig. \ref{FEx1}\emph{b}). Thus, $\cI=d=1$, $\vk_-=2$, $\vk_+=3$, $\n=4$, and the estimate \er{CEN111} has the
form \er{es1D}.\qq \BBox

\medskip

The estimates \er{esKL} and \er{es1D} are quite rude. Now we obtain
more accurate estimates for these graphs using Corollaries
\ref{CHEs} and \ref{CLEs}.

\begin{example}\lb{EKL1}
Let $\D$ be the combinatorial Laplacian defined by \er{DLO} on the
Kagome  lattice $\bK$ (see Fig.~\ref{FEx2}a). Then its total bandwidth $\gS(\D)$ satisfies
\[\lb{esKL1}
\textstyle2\leq\gS(\D)\leq12.
\]
\end{example}

\no \textbf{Proof.}  We estimate the total bandwidth for the
combinatorial  Laplacian $\D$ applying  the inequalities \er{CEH1}
as $n=1,2,3$. By \er{coCH}, we have
\[\lb{eqBB}
B_{n,1}\geq\cN_n^+, \qqq
B_{n,2}\geq2\cN_n^{odd},\qqq \forall\,n\in\N,
\]
where $\cN_n^+$ and $\cN_n^{odd}$ are defined in \er{cN+} and \er{cNo}.

There are no cycles of length one (i.e., loops) in the fundamental graph $\bK_*$ (see Fig.~\ref{FEx2}\emph{b}), i.e.,
\[\lb{Ncn1}
\cN_1^+=\cN_1^{odd}=0.
\]
Each oriented edge $\be\in\cA_*$ generates the cycle
$\bc_\be=(\be,\ul\be\,)$  of length 2 consisting of 1 backtrack
$(\be,\ul\be\,)$, but index of $\bc_\be$ is zero, since, by \er{cyin} and \er{inin},
$$
\t(\bc_\be)=\t(\be)+\t(\ul\be\,)=\t(\be)-\t(\be)=0.
$$
The graph $\bK_*$ also has the following proper cycles (i.e., cycles
without backtracking) of length 2:
$$
\bc_{2,1}=(\be_6,\be_3), \qqq \bc_{2,2}=(\be_1,\be_4), \qqq \bc_{2,3}=(\be_2,\be_5)
$$
with indices
$$
\t(\bc_{2,1})=(-1,0),\qqq \t(\bc_{2,2})=(0,1),\qqq \t(\bc_{2,3})=(1,-1),
$$
their cyclic permutations (two permutations for each cycle
$\bc_{2,s}$,  $s=1,2,3$) and their reverse cycles. Thus,
\[\label{Ncn2}
\cN_2^+=12,\qqq\cN_2^{odd}=8.
\]
There are no cycles of length 3 with backtracking in $\bK_*$, since
there are no loops and each backtracking contributes 2 in the cycle
length. Thus, all cycles of length 3 are proper cycles
$$
\begin{array}{llll}
\bc_{3,1}=(\be_1,\be_2,\be_3),\qq & \bc_{3,2}=(\be_1,\be_2,\ul\be_6),\qq &
\bc_{3,3}=(\be_2,\be_3,\ul\be_4),\qq & \bc_{3,4}=(\be_3,\be_1,\ul\be_5),
\\
\bc_{3,5}=(\be_6,\be_5,\be_4),  & \bc_{3,6}=(\be_6,\be_5,\ul\be_1),
& \bc_{3,7}=(\be_5,\be_4,\ul\be_3), & \bc_{3,8}=(\be_4,\be_6,\ul\be_2)
\end{array}
$$
with indices
$$
\begin{array}{llll}
\t(\bc_{3,1})=(0,0), \qq & \t(\bc_{3,2})=(1,0), \qq &
\t(\bc_{3,3})=(0,-1), \qq & \t(\bc_{3,4})=(-1,1),\\ \t(\bc_{3,5})=(0,0),
\qq & \t(\bc_{3,6})=(0,-1), \qq & \t(\bc_{3,7})=(1,0), \qq & \t(\bc_{3,8})=(-1,1),
\end{array}
$$
their cyclic permutations (three permutations for each cycle $\bc_{3,s}$,
 $s\in\N_8$) and their reverse cycles. Then,
\[\label{Ncn22}
\cN_3^+=36,\qqq\cN_3^{odd}=24.
\]
Using \er{eqBB} -- \er{Ncn22} and $v_+=0$, we obtain that the lower estimates in \er{CEH1} for $n=1,2,3$  have the form
$$
\gS(\D)\geq0, \qqq
\gS(\D)\geq2, \qqq \gS(\D)\geq1,
$$
respectively. For $n=1$ the lower estimate is trivial. For $n=2$ the
lower estimate  is better than for $n=3$. The upper estimate in
\er{esKL1} was proved in Example \ref{EKL}. \qq \BBox

\medskip

\begin{example}\lb{E1D1}
Let $\D_\gn$ be the normalized Laplacian defined by \er{DNLA} on the
periodic graph $\cG$ shown in Fig.~\ref{FEx1}a. Then its total bandwidth $\gS(\D_\gn)$ satisfies
\[\lb{es1D1}
\textstyle \frac49\leq\gS(\D_\gn)\leq\frac43\,.
\]
\end{example}

\no \textbf{Proof.} We estimate the total bandwidth of $\D_\gn$
applying  the inequalities \er{CEN1'} as $n=1,2,3,4$. There are no
cycles of length 1 (i.e., loops) in the fundamental graph  $\cG_*$ (see Fig. \ref{FEx1}\emph{b}).
Thus, \er{coC} gives
$$
B_{1,1}=B_{2,1}=0 \qq \textrm{as} \qq n=1.
$$
There are no proper cycles (i.e., cycles without backtracking) of
length 2 in $\cG_*$. Thus, all cycles of length 2 have the form
$\bc_\be=(\be,\ul\be\,)$ for some edge $\be\in\cA_*$ and
$\t(\bc_\be)=0$. Then  \er{coC} gives
$$
B_{2,1}=B_{2,2}=0 \qq \textrm{as}  \qq n=2.
$$
There are no cycles of length 3 with backtracking in $\cG_*$, since
there are no loops and each backtracking contributes 2 in the cycle
length. Thus, all cycles of length 3 are proper cycles
$$
\bc_1=(\be_1,\ul\be_3,\ul\be_2),\qqq \bc_2=(\be_1,\be_4,\be_5)
$$
with indices $\t(\bc_1)=\t(\bc_2)=1$ and weights
$\o_\gn(\bc_1)=\o_\gn(\bc_2)=\frac1{18}$\,, their cyclic
permutations (three permutations for each cycle $\bc_s$, $s=1,2$)
and their reverse cycles. Then, \er{coC} gives
$$
\textstyle B_{3,1}=\frac23\,,\qq B_{3,2}=\frac43\, \qq \textrm{as}  \qq n=3.
$$
At last, the graph $\cG_*$ has the following prime cycle of length 4:
$$
\bc=(\be_2,\be_3,\be_4,\be_5),
$$
and 8 prime cycles of the form $(\be,\be',\ul\be',\ul\be)$ for all
pairs  of edges $\be,\be'\in\cA_*$ such that the terminal vertex of
$\be$ coincides with the initial vertex of $\be'$ and $\be'\neq\ul\be$. All remaining cycles of length four are $2$-multiple of prime cycles of length
2. Each cycle of length 4 has zero index. Then,  by \er{coC} as
$n=4$, we obtain
$$
B_{4,1}=B_{4,2}=0.
$$
For $n=3$ the lower estimate in \er{CEN1'} has the form
$$
\textstyle\gS(\D_\gn)\geq \frac49\,.
$$
For $n=1,2,4$ the lower estimate in \er{CEN1'} is trivial.

Since the number $\gb$ defined in \er{gb} is equal to $\frac23$,
the upper estimate in \er{CEN1'} yields
$$
\textstyle\gS(\D_\gn)\leq 2\cdot\frac23=\frac43\,. \qqq \BBox
$$

\begin{remark}
Comparing the estimates \er{esKL} and \er{es1D} with the estimates
\er{esKL1} and  \er{es1D1}, we see that the last ones are better.
Moreover, the upper estimate in \er{es1D1} is sharp, see remark
after Example \ref{E1D}.
\end{remark}

\medskip

\no\textbf{Acknowledgments.} \footnotesize Our study was supported
by the   RFBR grant No. 19-01-00094.

\end{document}